\documentclass[12pt,a4paper]{article}
\ifcase0
\or\documentclass[12pt,a4paper]{article}
\or\documentclass[12pt,a5paper]{article}
\fi

\usepackage{euler,amsmath,amsfonts,defns,url,graphicx}
\IfFileExists{ajr.sty}{\usepackage{ajr}}{}
\allowdisplaybreaks

\newcommand{\cL}{\mathcal{L}}
\newcommand{\cG}{\mathcal{G}}

\newcommand{\spde}{\textsc{spde}}
\newcommand{\ibc}{\textsc{cc}}
\newcommand{\scm}{\textsc{scm}}
\newcommand{\ssm}{\textsc{ssm}}

\newcommand{\inti}{\int_{-\infty}}

\newcommand{\up}{U_{j+1}}
\newcommand{\um}{U_{j-1}}

\newcommand{\uj}{U_{j}}
\newcommand{\duj}{\D t{U_{j}}}
\newcommand{\xj}{X_{j}}
\newcommand{\xpm}{X_{j\pm 1}}
\newcommand{\pars}{\vec\epsilon} 
\newcommand{\cosin}{\operatorname{csn}}
\newcommand{\Z}[1]{\operatorname{\mathcal{H}}_{#1}}

\newcommand{\mud}{\mu\delta}

\newcommand{\bspde}{Burgers' \spde~\eqref{eq:burg}}
\newcommand{\spdee}{\spde~\eqref{eq:spde}}
\newcommand{\eigf}{v} 
\newcommand{\vecs}[1]{\mathbb{#1}}
\def\arcproject#1{}\def\arcstatement#1{}

\title{Resolve subgrid microscale interactions to discretise stochastic partial differential equations}

\author{A.~J. Roberts\thanks{School of Mathematical Sciences, University of Adelaide, South Australia~5005, Australia.  
\protect\url{mailto:anthony.roberts@adelaide.edu.au}}}

\begin{document}

\maketitle

\begin{abstract}
Constructing discrete models of stochastic partial differential equations is very delicate. Stochastic centre manifold theory provides novel support for coarse grained, macroscale, spatial discretisations of nonlinear stochastic partial differential  or difference equations such as the example of the stochastically forced Burgers' equation. Dividing the physical domain into finite length overlapping elements empowers the approach to resolve fully coupled dynamical interactions between neighbouring elements. The crucial aspect of this approach is that the underlying theory organises the resolution of the vast multitude of subgrid microscale noise processes interacting via the nonlinear dynamics within and between neighbouring elements. Noise processes with coarse structure across a finite element are the most significant noises for the discrete model.  Their influence also diffuses away to weakly correlate the noise in the spatial discretisation. Nonlinear interactions have two further consequences: additive forcing generates multiplicative noise in the discretisation; and effectively new noise processes appear in the macroscale discretisation. The techniques and theory developed here may be applied to soundly discretise many dissipative stochastic partial differential and difference equations.
\end{abstract}


\tableofcontents

\section{Introduction}

This article develops a systematic approach to constructing spatially discrete models of stochastic partial differential equations and stochastic semi-discrete equations (collectively denoted \spde{}s).  Although this article considers differential examples of the theory, following analogous deterministic results~\cite{Roberts08c}, section~\ref{sec:cm} comments how the theory also underpins the modelling on a coarse grid of stochastic dynamics on a fine grid, and thence potentially across a multigrid hierarchy.  Thus the research here contributes to understanding how stochastic information is transferred over levels of scale---an outstanding need as reported by Dolbow et al.~\cite{Dolbow04} who summarised a wide ranging workshop.  More recently, reporting a review of Applied Mathematics, Brown et al.~\cite[\S2.1.3]{Brown08} similarly called  for ``new approaches for efficient modeling of large stochastic systems''.  This article explores sound modelling of complex stochastic systems over large space domains. 

The aim is to use stochastic centre manifold (\scm) theory and techniques to ensure the accuracy, stability and efficiency of relatively coarse spatial discretisations of \spde{}s.  Understanding such \scm\ theory is best obtained via stochastic normal form theory~\cite{Arnold03} developed because of the success of deterministic normal forms~\cite[e.g.]{Murdock03}:  a stochastic coordinate transform separates the coarse grid slow  dynamics in the system from the fast subgrid dynamics; in the new coordinates the subgrid modes decay exponentially quickly from all initial conditions in a finite domain;  leaving the slow modes to evolve on the stochastic slow manifold as a rigorous closure on the coarse grid.  Related dynamical systems approaches by Blomker, Hairer and Pavliotis~\cite{Blomker04} supports stochastic Ginzburg--Landau model in pattern forming systems.  Caraballo, Langa and Robinson~\cite{Caraballo01} and Duan, Lu \& Schmalfuss~\cite{Duan04} proved the existence of invariant manifolds for multiplicative stochastic dynamics.
Here, due to the stochastic effects over many length and time scales, the practical challenge is to close on a macroscale grid the intricate spatio-temporal dynamics of a \spde. Numerical methods to integrate stochastic \emph{ordinary} differential equations are known to be delicate and subtle~\cite[e.g.]{Kloeden92}. We surely need to take considerable care for \spde{}s as well~\cite[e.g.]{Grecksch96, Werner97} in order to model the subgrid microscale nonlinear interactions.  Modern dynamical systems theory provides a requisite careful methodology. 

The sound methodology for modelling \spde{}s presented here has the potential to find future applications underpining multiscale modelling in areas such as climate, materials science and the biosciences~\cite[e.g.]{Dolbow04}. For example, the gap-tooth scheme of Kevrekidis et al.~\cite{Gear03, Samaey03a, Samaey03b} is often implemented with particle simulators. Coupling patches of the stochastic microscale dynamics of such particle simulators, analogous to the coupling of stochastic element developed here, will be an important generalisation of coupling patches with deterministic dynamics~\cite{Roberts04d} to the stochastic dynamics explored herein.

Earlier research developed a dynamical systems approach to discretising deterministic \pde{}s~\cite[e.g.]{MacKenzie05a} and briefly looked at the specific case of stochastic \emph{linear} diffusion~\cite{Roberts06g}.  Others proved the appearance of stochastic Ginzburg--Landau models in pattern forming systems~\cite{Blomker04}, and the existence and uniqueness of global solutions to the \bspde~\cite{DaPrato94, DaPrato96}. Section~\ref{sec:cm} significantly extends new theoretical support to the dynamical systems discrete modelling of the general \emph{nonlinear} \spde
\begin{equation}
    \D tu=\cL(u) u+\alpha f(u)+\sigma\phi(u,x,t) \,,
    \label{eq:spde}
\end{equation}
(adopting the notation of most physicists and engineers) for a stochastic field~$u(x,t)$ evolving in time~$t$ in one spatial dimension, where $\cL(u)$~is a smooth, second order, dissipative operator, $f(u)$~denotes smooth, deterministic, perturbations (possibly involving derivatives of the field~$u$), and $\phi$~denotes the stochastic effects in the \spde.  When necessary for definite theoretical statements and for displayed numerical simulations, we adopt boundary conditions for the \spdee\  of $L$-periodicity, $u(x+L,t)=u(x,t)$, and assume initial conditions lie within the finite domain of validity of the \scm\ theory.  Section~\ref{sec:cm} discusses that the \spdee\ not only includes the class of \spde{}s differential in space, but also, such as when $\cL$~and~$f$ are discrete operators on a lattice, a wide class of microscale stochastic  discrete difference equations. The analysis also applies to many microscale lattice stochastic rules.

\begin{figure}
    \centering
    \includegraphics{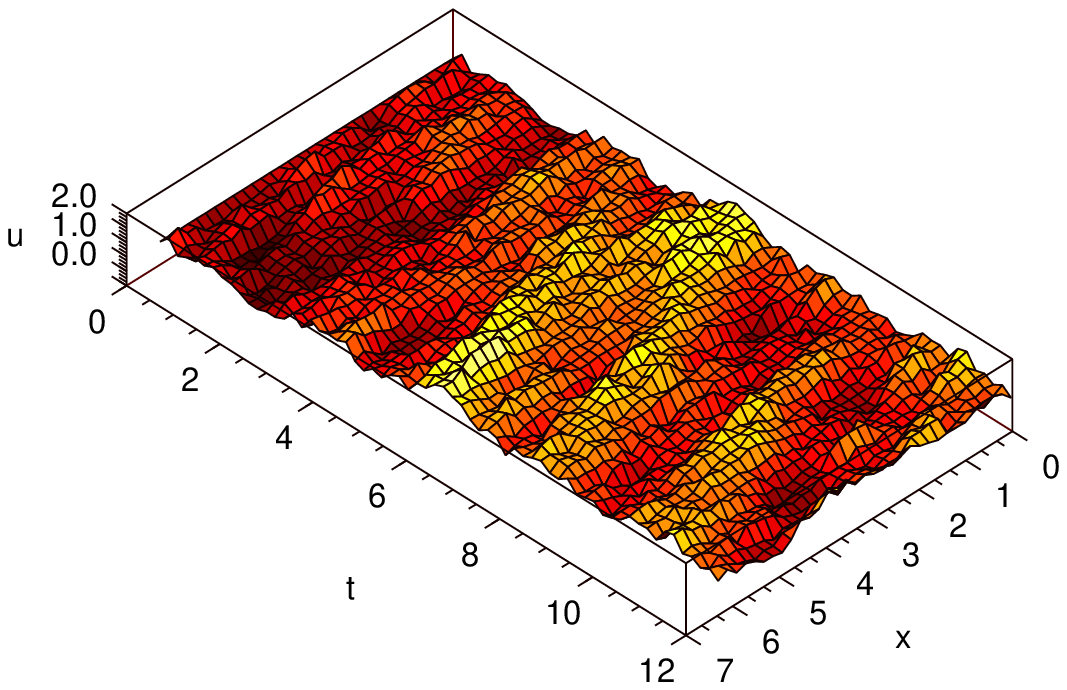}
	\caption{an example microscale simulation of one realisation,~$u(x,t)$, of the stochastically forced \bspde, with nonlinearity $\alpha=3$ and noise intensity $\sigma=1$\,.  This simulation uses a fine space-time mesh with $\Delta x=\pi/16$ and $\Delta t=0.01$ but plotted every $19$th~time step.}
    \label{fig:micro}
\end{figure}

Typical issues and results are illustrated herein via the definite example of discretising in space the stochastically forced \emph{nonlinear} Burgers' equation:
\begin{equation}
    \D tu +\alpha u\D xu =\DD  x u +\sigma\phi(x,t)
    \label{eq:burg}\,.
\end{equation}
The simulation of the \bspde\ in Figure~\ref{fig:micro} illustrates complicated microscale fluctuations and their cumulative appearance in the macroscale. Givon et al.~\cite[p.R58]{Givon04} similarly used specific example problems to develop and illustrate many issues in their approach. The aim here is to establish stochastic centre manifold (\scm) theory and techniques as a practical approach to ensure the accuracy, stability and efficiency of numerical discretisations of \spde{}s. 

Even discretising linear \spde{}s has subtleties. Consider stochastically forced, linear diffusion in space (\bspde\ with $\alpha=0$). The simplest finite difference approximation on a regular grid in~$x$, say $\xj =jh$ for some constant grid spacing~$h$, is
\begin{equation}
	\duj = \frac{\up-2\uj+\um}{h^2}+\sigma\phi(\xj,t)\,,
    \label{eq:diffmm}
\end{equation}
where $\uj$ is the value of the field~$u(x,t)$ at the grid points~$\xj$.  However, the systematic analysis of specific earlier research~\cite{Roberts06g} recommends we use instead
\begin{equation}
	\duj \approx \frac{\up-2\uj+\um}{h^2}
        +\sigma\left[ \sqrt{\rat57}\psi_j
        -\rat1{24}\sqrt{\rat75}(\psi_{j+1}-2\psi_j+\psi_{j-1}) 
        +\sqrt{\rat27}\hat\psi_j
        \right],
        \label{eq:diffm}
\end{equation}
for some independent noise processes $\psi_j$~and~$\hat\psi_j$.  That is, the point sample~$\phi_j(t)=\phi(\xj,t)$ of the noise, as invoked in~\eqref{eq:diffmm}, is better replaced by two components: a component~$\hat\psi_j(t)$ which is uncorrelated across the grid points; and a component~$\psi_j(t)$ whose direct influence distributes over the evolution of neighbouring grid values. The rationale is that spatial diffusion on the subgrid microscale between the macroscale grid points weakly correlates the noise applied to each grid value~\cite{Roberts06g}. We use the stochastic centre manifold theory of Boxler~\cite{Boxler89} to account explicitly and systematically for such subgrid microscale physical processes and interactions. 

Subgrid microscale processes are especially important for nonlinear systems.  One example is furnished by the stochastic, Burgers'-like, spatially discrete, fine grid equation
\begin{equation}
\frac{du_i}{dt}=\frac{4}{h^2}(u_{i+1}-2u_i+u_{i-1})
-\frac\alpha hu_i(u_{i+1}-u_{i-1}) +\sigma\phi_i(t) ,
\label{eq:bssde}
\end{equation}
where subscript~$i$ indexes the fine grid, the first term on the right-hand side is a diffusive-like dissipative reaction with neighbours, the second term represents a form of nonlinear advection of strength~$\alpha$, and the last term is an additive stochastic forcing independent at each grid point but of uniform strength~$\sigma$.  The theory described in Section~\ref{sec:cm} for the \spdee\ also supports the modelling of~\eqref{eq:bssde} on a coarser grid with twice the spacing.  Extending earlier research for deterministic systems~\cite{Roberts08c}, define coarse grained amplitudes~$\uj(t)=u_{2j}(t)$, our coarser grid model of~\eqref{eq:bssde} is
\begin{equation}
\duj\approx \frac{1}{h^2}(\up-2\uj+\um)
-\frac\alpha {2h}\uj(\up-\um) +\sigma\left[\psi_{j0}(t)
-\frac{\alpha h}8\uj\psi_{j1}(t)\right].
\label{eq:bssdm}
\end{equation}
To this level of approximation the coarse grid evolution~\eqref{eq:bssdm} is  the same form as the fine grid~\eqref{eq:bssde} with appropriately renormalised diffusion and nonlinear advection, but with noise now weakly correlated with its neighbours as the coarse grid noise is the multigrid restriction $\psi_{j0}=\frac14\phi_{2j-1}+\frac12\phi_{2j}+\frac14\phi_{2j+1}$\,, and with a new multiplicative noise arising from resolving the subgrid interaction between the nonlinear advection and structure in the subgrid noise.  In contrast to other methods, stochastic centre manifold theory~\cite{Boxler89} provides a systematic approach to resolving the macroscale influence of such microscale stochastic nonlinear interactions.

Section~\ref{sec:sto} provides one resolution of the challenge to model mean effects on the grid values~$\uj$, the macroscale, induced by the microscale processes and interactions, both linear and nonlinear.  This issue of macroscale, grid value modelling closure is a longstanding problem in even deterministic dynamics, let alone a reportedly outstanding need in stochastic dynamics.  For example,  should the advection term~$\alpha uu_x$ in \bspde\ be discretised directly? or as $\rat12\alpha(u^2)_x$? or as the $1:2$~combination of these suggested by Fornberg~\cite{Fornberg73} and used to improve stability~\cite[e.g.]{Foias91}. Which is better? The traditional approach of considering the discrete modelling of each term separately does not tell us. Instead, to find a sound discretisation we consider the interaction of all terms in the \spde\ by resolving subgrid scale processes and hence determine how the macroscale grid values can capture these subgrid scale processes.  Section~\ref{sec:onhnni} discusses a case of \bspde\ where traditional straightforward numerical approximations miss all of the easily apparent effects of the microscale noise on the macroscale discretisation.   Here we develop further and validate the stochastic centre manifold (\scm) methodology for closure of macroscale discretisations. The wide ranging reports of Dolbow et al.~\cite{Dolbow04} and Brown et al.~\cite{Brown08} identify such closure as an outstanding challenge in multiscale modelling of stochastic physical systems of great interest to applied mathematics.

\subsection{Divide space into overlapping finite elements}
\label{sec:dsiofe}

The method of lines discretises a deterministic \pde\ in space~$x$ and integrates in time as a set of ordinary differential equations---sometimes called a semi-discrete scheme \cite[e.g.]{Foias91, Foias91b}.  Similarly, this article only discusses the relatively coarse grained spatial discretisation of \spde{}s and treats the resulting set of stochastic ordinary differential equations, such as the \sde~\eqref{eq:diffm}, as a continuous time, stochastic dynamical system.  Earlier research~\cite{Roberts05c} explored how to model the \bspde\ in just a small domain of one element and focussed on how to handle the complicated microscale interactions of noise in an \spde.  This article not only establishes theoretical support for a wider range of \spde{}s, it focusses on the vastly more complicated nonlinear interactions of noise not only within but also between near neighbouring spatial elements.

Place the spatial discretisation of a nonlinear \spdee, such as \bspde, within the purview of \scm\ theory by the following artifice. Let equi-spaced grid points at~$\xj$ be a distance~$h$ apart. The $j$th~element is notionally $|x-\xj|<h/2$ but we extend and overlap the elements to be $|x-\xj|\leq h$\,.  Let there be $m$~such elements. Following the analogous and proven  approach to discretising deterministic \pde{}s~\cite[e.g.]{MacKenzie05a, Roberts00a, MacKenzie09b}, form \emph{overlapping} finite elements by introducing artificial, nonlocal, coupling conditions (\ibc{}s) between neighbouring elements:
\begin{equation}
    u_j(\xpm,t) = (1-\gamma)u_j(\xj,t) + \gamma u_{j\pm1}(\xpm,t) ,
    \label{eq:ibc}
\end{equation}
where $u_j(x,t)$ denotes the subgrid microscale field of the $j$th~element. The coupling parameter~$\gamma$ controls the flow of information between adjacent elements: when $\gamma=0$\,, adjacent elements are decoupled by the \ibc~\eqref{eq:ibc}; whereas when $\gamma=1$\,, the \ibc~\eqref{eq:ibc} requires the field in the $j$th~element to extrapolate to the neighbouring elements' field at their grid point.  Others found similarly overlapping elements useful in several multiscale methods~\cite[e.g.]{E04, Samaey03b, Gander98}. The specific coupling \ibc{}s~\eqref{eq:ibc} used here ensure discrete models are consistent with deterministic \pde{}s to high order in~$h$ as the element size~$h\to0$\,, both linearly~\cite{Roberts00a}, and nonlinearly and in multiple space dimensions~\cite[\S4]{MacKenzie09b}.

Importantly, Section~\ref{sec:cm} establishes that the methodology developed with the \ibc{}s~\eqref{eq:ibc} applies to a wide variety of continuum and discrete \spde{}s, not just the example \bspde.

\subsection{Model nonlinear stochastic dynamics}

This paragraph focuses on nonlinear interactions rather than stochastic interactions between elements. Sections \ref{sec:quad}~and~\ref{sec:sto} construct the following low order, discrete system of \sde{}s of the nonlinear dynamics of the stochastically forced \bspde:
\begin{eqnarray}
    \duj&\approx&
    \frac1{h^2}\big( 1+\rat1{12}\alpha^2h^2\uj^2 \big)(\up-2\uj+\um)
    -\alpha\frac{1}{2h}\uj(\up-\um)
    \nonumber\\&&{}
    +\sigma\left[ \phi_{j,0} -\alpha\frac{2h}{\pi^2}\phi_{j,1}\uj
    -\alpha^2\frac{8h^2}{3\pi^4}\phi_{j,2}\uj^2 \right]
    + .01643\alpha^2h^2\sigma^2\uj\,.
    \label{eq:lowg}
\end{eqnarray}
This model is for the case when the subgrid microscale noise within each element is truncated to the first three Fourier modes:
\begin{displaymath}
    \phi(x,t)=\phi_{j,0}(t) +\phi_{j,1}(t)\sin[\pi(x-\xj)/h]
    +\phi_{j,2}(t)\cos[2\pi(x-\xj)/h]\,.
\end{displaymath}
The first line of the discretisation~\eqref{eq:lowg} is the so-called holistic discretisation for the deterministic Burgers' equation which has good properties on finite sized elements~\cite{Roberts98a}; in particular, the nonlinearly enhanced diffusion promotes stability of the scheme for non-small field~$u$. The second line of the discretisation~\eqref{eq:lowg} approximates some of the stochastic influences. Analogous to the coarse graining modelling \eqref{eq:bssde}$\mapsto$\eqref{eq:bssdm}, the nonlinearity in the subgrid microscale dynamics of Burgers' equation transforms the additive noise of the \spde~\eqref{eq:burg} into multiplicative noise components in the discretisation, such as $\phi_{j,1}\uj$.  Many modelling schemes miss such multiplicative noise terms because they do not resolve the subgrid microscale processes that are revealed by the \scm\ methodology.

As seen in Figure~\ref{fig:micro}, stochastic forcing typically generates microscale spatial structures with high wavenumber, steep variations. Through a form of stochastic resonance, the life-time of these modes may be important on the large scale dynamics. For example, the last term of the discrete model \sde~\eqref{eq:lowg}, being proportional to~$\sigma^2\uj$, arises from nonlinear self-interactions of subgrid scale noise, as discussed in Section~\ref{sec:stow}. Herein the term `stochastic resonance' includes phenomena where stochastic fluctuations interact with each other and themselves through nonlinearity in the dynamical system to generate not only long time drifts but also potentially to change stability~\cite[e.g.]{Knobloch83, Boxler89, Drolet01, Roberts03c, VandenEijnden05c}. Coarse grained discrete model \sde{}s, such as \eqref{eq:lowg}~and~\eqref{eq:diffm}, which use large space-time grids for efficiency, must account for such subgrid microscale dynamics in their closure in order to resolve the significant subtle interactions.

Stochastic averaging similarly derives mean effects on slow modes from fast time stochastic effects~\cite[e.g.]{Pavliotis07}, and complicated deviation arguments estimate fluctuations~\cite{Wang2008}.  However, in discretising the \spdee\ there is only a finite scale separation between slow and fast modes: in forming discrete models there is no asymptotically small parameter~$\epsilon$ measuring an asymptotically infinite scale separation as required by stochastic singular perturbation theory.  Dolbow et al.~\cite[p.30]{Dolbow04} specifically call for ``new multiscale mathematical methods developed and used to derive multiscale models for some of the `difficult' cases in multiscale science; e.g., problems without strong scale separation,''.  Our application of \scm\ provides such a methodology.

Stable implicit integration schemes do not resolve subgrid microscale 
fluctuations.  For example, recall that the well
established Crank--Nicholson scheme for spatio-temporal dynamics is
based upon the stability of the implicit scheme 
\begin{displaymath}
    \frac{u_{t+\Delta t}-u_t}{\Delta t}=-\beta\frac{u_{t+\Delta t}+u_t}2
    \qtq{for the \ode}
    \dot u=-\beta u\,.
\end{displaymath}
The solution is indeed stable,
\begin{displaymath}
	u(t)=u(0)\left( \frac{1-\beta\Delta t/2} {1+\beta \Delta t/2}
	\right)^{t/\Delta t}\,,
\end{displaymath}
and accurate for small~$\beta\Delta t$.
Thus in general such schemes preserve noise characteristics,
such as drift and volatility, of the evolution of \emph{slow}
macroscopic modes which have small rates~$\beta$.
However, the dynamics of microscale modes, with large~$\beta\Delta t$,
are badly misrepresented by such schemes when employing a macroscopic
time step of relatively large~$\Delta t$; for example, any exponential
decay of microscale modes is badly simulated by near oscillations in the Crank--Nicholson scheme.  Such bad representation may be acceptable in deterministic \pde{}s where the microscale modes are not forced.  But in a \spde, with broad spectrum noise, the microscale modes are continually forced and their life-time requires reasonable resolution. 
The methodology proposed here provides a systematic method
to derive macroscale closure of  the microscale~noise.

Stochastic centre manifold theory supports the large time macroscopic modelling of detailed stochastic microscopic dynamics, even when there is a finite scale separation between `slow' and `fast' modes. Knobloch \& Wiesenfeld~\cite{Knobloch83} and Boxler~\cite{Boxler89, Boxler91}, for example, explicitly used \scm\ theory to support the modelling of \sde{}s and \spde{}s. Boxler~\cite{Boxler89} proved that ``stochastic center manifolds, share all the nice properties of their deterministic counterparts''. Many, such as Berglund \& Gentz~\cite{Berglund03}, Bl\"omker, Hairer \& Pavliotis~\cite{Blomker04} and Kabanov \& Pergamenshchikov~\cite{Kabanov03}, use the same separation of time scales that underlies \scm\ theory to form and support low~dimensional, long time models of \sde{}s and \spde{}s that have both fast and slow modes. Similarly, Bensoussan \& Flandoli~\cite{Bensoussan95} proved the existence and relevance of finite dimensional stochastic inertial manifolds for a wide class of stochastic systems in a Hilbert space.  Stochastic normal form theory of Arnold~\cite[e.g.]{Arnold03, Roberts06k} also supports the existence and relevance of stochastic centre manifolds.  Centre manifold theory also supports the discretisation of deterministic partial differential equations~\cite{Roberts98a, Roberts00a, MacKenzie00a, Roberts01a, Roberts01b, MacKenzie03}. Section~\ref{sec:cm} merges these two applications of theory to model \spde{}s with sound theoretical support. Coupling neighbouring finite elements together forms macroscale discrete models of \spde{}s such as \eqref{eq:diffm}, \eqref{eq:bssdm} and~\eqref{eq:lowg}. The major complication is to account for noise and its subgrid dynamics which are distributed independently across space as well as time, both within elements and between neighbouring elements. Computer algebra---given in full detail, documented and freely available for checking and modification on an institutional preprint server~\cite{Roberts06b}---conveniently handles the fearsome details of the nonlinear subgrid dynamics and the interelement interactions in the stochastic centre manifold modelling.

We discuss the forcing~$\phi(x,t)$ as a white noise that is delta correlated in both space and time. However, computational limitations often require the truncation to a few Fourier modes as in~\eqref{eq:lowg}.  Nonetheless, most of the analysis and models in Sections \ref{sec:cm}~and~\ref{sec:quad} also hold for deterministic forcing~$\phi(x,t)$.  Further, correlated noise~$\phi$ could be handled in the same dynamical systems approach by introducing auxiliary \spde{}s such as
\begin{displaymath}
    \D t\phi=-\cG\star\phi+\xi(x,t)\,,
\end{displaymath}
for some white space-time noise~$\xi$ and convolution kernel~$\cG(x)$. The stochastic~$\phi$ would then have spatial correlations with power spectrum${}\propto 1/|\tilde\cG(k)|^2$ for Fourier transform~$\tilde\cG$ of the kernel~$\cG(x)$.  Consequently, this approach encompasses quite general noise processes.

Herein, interpret all noise processes and equations in the Stratonovich sense so that the rules of traditional calculus apply. Thus the direct application of this modelling is to physical and engineering systems where the Stratonovich interpretation is the norm.

\section{Stochastic centre manifold theory underpins modelling}
\label{sec:cm}

This section details one way to place the spatial discretisation of \spde{}s within the purview of stochastic centre manifold (\scm) theory.  The theory was elaborated by Arnold~\cite[Chapt.~7]{Arnold03}, or see the freely available concise summary~\cite[Appendix~A, e.g.]{Roberts05c}. Subject to some conditions, extant \scm\ theory assures us of the existence and relevance of discrete models of the general \spdee.  When a definite example is discussed we invoke the \bspde, or the discrete microscale dynamics of~\eqref{eq:bssde}.  Earlier related work~\cite{Roberts03c} was limited to linear differential dynamics and invoked different, less practical, interelement coupling conditions.  To proceed, first we embed the dynamics of the \spdee\ into a higher dimensional stochastic system. Then from the base of a subspace of equilibria, this section establishes the existence and emergence of a stochastic slow manifold that forms the macroscale discrete model.  The key here is simply to establish the preconditions for the application of stochastic centre manifold theory.  Then the established \scm\ theory rigorously supports a wide range of models and illuminates a wide range of modelling issues---to obtain useful results we only need to verify the preconditions.

The first step is to embed the dynamics of the \spdee\ onto the overlapping elements of Section~\ref{sec:dsiofe}.  Recall that $u_j(x,t)$~denotes the stochastic field in the $j$th~element $X_{j-1}\leq x\leq X_{j+1}$\,.  Then on each element solve the \spdee\ with the field in the $j$th~element, namely we analyse the set of \spde{}s
\begin{equation}
    \D t{u_j}=\cL(u_j) u_j+\alpha f(u_j)+\sigma\phi_j(u_j,x,t) \,,
    \label{eq:spdei}
\end{equation}
where $\phi_j$ is the noise term on the $j$th~element.  Couple the fields in neighbouring elements with the \ibc~\eqref{eq:ibc}.  Such coupling of overlapping elements appears analogous to the `border regions' of the heterogeneous multiscale method~\cite[e.g.]{E04}, to the `buffers' of the gap-tooth scheme~\cite[e.g.]{Samaey03b}, and to the overlapping domain decomposition that improves convergence in waveform relaxation of parabolic \textsc{pde}s~\cite[e.g.]{Gander98}.  When necessary for definitive theory and for numerical simulations, the $L$-periodic boundary conditions ($L=mh$) for the global field~$u(x,t)$ require that the fields in the $m$~elements satisfy $u_{j\pm m}(x\pm L,t)=u_j(x,t)$ for all~$x$ and elements~$j$.   Equivalently, consider the \spde\ in the interior of a domain large enough so the physical boundaries are far enough away to be immaterial.  That is, the analysis throughout assumes space and time are statistically homogeneous so that the resultant discretisations are statistically homogeneous as seen in \eqref{eq:diffm}, \eqref{eq:bssdm} and~\eqref{eq:lowg}.   Evidently, the set of \spde{}s~\eqref{eq:spdei} form a higher dimensional stochastic system that effectively reduces to \spdee\ in the limit of full coupling.

I avoid defining a precise probability space for analysis of the \spde{}s~\eqref{eq:spdei}.  The reason is that extant theorems place differing conditions of the nature of the functions $\cL(u)$, $f(u)$ and~$\phi(u,x,t)$ appearing in the \spde; see Corollary~\ref{lem:exist} for examples.  Herein we primarily develop a formal methodology applicable to all established rigorous theoretical conditions, but flexible enough to cater also for a wide class of physically interesting stochastic systems.     

Second, we anchor the discrete modelling upon the subspace of equilibria~$\vecs E_0$ of piecewise constant solutions, $u_j=\uj ={}$constant, with no noise, $\sigma=0$\,, no interelement coupling, $\gamma=0$ in the \ibc~\eqref{eq:ibc}, and no nonlinearity, $\alpha=0$\,. For each of the equilibria in~$\vecs E_0$ it is straightforward to find the linear Oseledec spaces that form the foundation of the \scm~\cite[p.212]{Roberts05c} as here they are then standard linear eigenspaces. Linearising about each equilibria, for small perturbations~$u'_j(x,t)$ the dynamics of the \spdee\ with coupling conditions~\eqref{eq:ibc} reduce to that of dissipation within each element isolated from its neighbours:
\begin{equation}
    \D t{u'_j}=\cL_j u'_j  
    \quad\mbox{such that}\quad
    u'_j(\xpm,t)=u'_j(\xj,t),
    \label{eq:spdel}    
\end{equation}
where the prime denotes a perturbation to the equilibrium field of each element of~$\vecs E_0$ and $\cL_j=\cL(\uj )$;
for example, \bspde\ linearises to the diffusion equation $u'_t=u'_{xx}$\,. 

The linearised problem~\eqref{eq:spdel} has dissipative dynamics on each of the $m$~elements.
\emph{Assume the spectrum~$\{0, -\beta_{j1}, -\beta_{j2}, \ldots\}$ of the dissipative operator~$\cL_j$ is discrete with negative real parts bounded away from zero (apart from the single zero eigenvalue):}
$0>-\beta\geq\Re(-\beta_{j1})>\Re(-\beta_{j2})>\cdots$ for all elements~$j$. For example, the linearised problem~\eqref{eq:spdel} for \bspde\ is spatial diffusion which within each element has spectrum given by the negative of the decay rates $\beta_k={\pi^2k^2}/{h^2}$\,, for integer $k=0,1,2,\ldots$\,; the corresponding (generalised) eigenfunctions in each element are the modes\footnote{The corresponding adjoint (generalised) eigenfunctions are $(1-|\theta|/\pi)\cos k\theta$\,, $\sin k\theta$ and $\sin k|\theta|$ because the adjoint boundary conditions to those in~\eqref{eq:spdel} are $u'_j=0$ at $x=\xpm$\,, $u'_j$~is continuous at~$\xj$, and the nonlocal derivative condition $u'_{jx}(\xj^+,t)-u'_{jx}(\xj^-,t)=u'_{jx}(X_{j+1},t)-u'_{jx}(X_{j-1},t)$.}
\begin{equation}
\cos k\theta\text{ \ for even }k,\quad
\sin k\theta\text{ \ for }k\geq1\,,\quad
\theta\sin k\theta\text{ \ for even }k\geq2\,,
\label{eq:eigfns}
\end{equation}
where the subgrid variable $\theta=\pi(x-\xj)/h$~measures subgrid position relative to the centre grid point within each element (the $j$th~element lies between $\theta=\pm\pi$).
The $k=0$ mode, $u\propto{}$constant in each element, is linearly neutral as its decay rate~$\beta_0=0$\,. Thus, linearised about each equilibria in the subspace~$\mathbb E_0$, subgrid structures within each element decay so that a global piecewise constant field emerges exponentially quickly---at least as fast as~$\exp({-\beta t})$.

\scm\ theory asserts that such emergence is robust to nonlinear and stochastic perturbations.  To address stochastic effects we need to expand the noise in a complete set of basis functions.  Because we embed the dynamics in overlapping elements with `overlapping' eigenfunctions, such as~\eqref{eq:eigfns}, choose a subset~$\eigf_k(x-\xj)$ of the eigenfunctions that are complete over the non-overlapping domains $|x-\xj|<h/2$; for \bspde\ decompose the noise within each element as a linear combination of 
\begin{equation}
    \eigf_k(x-\xj)= \cosin k\theta=\left\{
    \begin{array}{ll}
        \cos k\theta\,, &\text{for even }k,\\
        \sin k\theta\,, &\text{for odd }k.
    \end{array}
    \right.
    \label{eq:four}
\end{equation}
This decomposition is analogous to Example~5.2.2 of Da~Prato \& Zabczyk~\cite[see also p.259]{DaPrato96}. Thus additive noise in the element \spde~\eqref{eq:spdei} for the \bspde\ is
\begin{equation}
    \phi_j(x,t)=\sum_{k=0}^\infty \phi_{j,k}(t)\eigf_k(x-\xj)
    =\sum_{k=0}^\infty \phi_{j,k}(t)\cosin k\theta\,,
    \label{eq:onoise}
\end{equation}
where $\phi_{j,k}$ denotes the noise process of the $k$th~wavenumber in the $j$th~element. \emph{Assume the noise processes~$\{\phi_{j,k}\}$ are independent.} Simple numerical methods, such as Galerkin projection onto the coarsest mode~$\eigf_0$, would ignore the `high wavenumber' subgrid modes~$\eigf_k$, $k\geq1$\,, of the noise~\eqref{eq:onoise} and hence miss subtle but potentially important subgrid and inter-element interactions such as those seen in the models \eqref{eq:diffm}, \eqref{eq:bssdm} and~\eqref{eq:lowg}. In contrast, the systematic nature of our application of \scm\ theory accounts for subgrid microscale interactions as an asymptotic series in the noise amplitude~$\sigma$, inter-element coupling~$\gamma$ and  nonlinearity~$\alpha$.

\subsection{A stochastic slow manifold exists} 

The nonlinear forced \spde~\eqref{eq:spdei} with inter-element coupling conditions~\eqref{eq:ibc} linearises to the dissipative \pde~\eqref{eq:spdel}. To account for non-zero parameters, adjoin the system of three trivial \textsc{de}s $d\pars/dt =\vec 0$\,, where $\pars=(\sigma,\gamma,\alpha)$\,. In the extended state space~$(\vec u,\pars)=(u_1,\ldots,u_m,\sigma,\gamma,\alpha)$, the linearised \pde\ has $m+3$~eigenvalues of zero and all other eigenvalues have negative real part~$\leq-\beta$: this upper bound is $-\pi^2/h^2$ for the example of \bspde. Thus the \spde~\eqref{eq:spdei} has an $m+3$~dimensional slow subspace. Because of the pattern of the eigenvalues, \scm\ theory~\cite[Theorem~3, p.212]{Roberts05c} assures us a corresponding $m+3$~dimensional stochastic slow manifold (\ssm) exists under certain conditions.

\begin{corollary}[existence]\label{lem:exist}
For Lipschitz nonlinearities~$\cL(u)$,  $f(u)$~and~$\phi(u)$, and in some finite neighbourhood of the equilibria,~$(\mathbb E_0,\vec 0)$, there exists an $m+3$~dimensional \ssm\ for the general \spde~\eqref{eq:spdei} in which the field in the $j$th~element is $u_j=w_j(\vec U(t),x,t,\pars)$ where the $j$th component~$\uj$ of vector~$\vec U$ measures the amplitude of the neutral mode~$\eigf_0(x-\xj)$ in the $j$th~element, and where the amplitudes~$\uj$ evolve according to $\dot \uj =g_j(\vec U,t,\pars)$ for some function~$g_j$, but provided
\begin{itemize}
	\item either the nonlinear \spde~\eqref{eq:spdei} is effectively finite dimensional\footnote{A \spde\ is effectively finite dimensional if there exists a wavenumber~$K$ such that modes~$\eigf_k$ for $k\geq K$ do not affect, through~$\cL$, $f$~or~$\phi$, the dynamics for modes with $0\leq k<K$\,.}~\cite[Theorems 5.1~and~6.1]{Boxler89}; or
    

    \item  the noise in the \spde~\eqref{eq:spdei} is multiplicatively linear in~$u$,
    $\phi=u\psi(x,t)$ \cite[Theorem~A]{Wang06}.

\end{itemize}
\end{corollary}
Analogously, Bl\"omker et al.~\cite[Theorem~1.2]{Blomker04} rigorously proved the existence and relevance of a stochastic Ginzburg--Landau model to the `infinite dimensional' stochastic forced Swift--Hohenberg
\pde. Similarly, Caraballo, et al.~\cite{Caraballo01} and Duan et al.~\cite{Duan04} analysed the existence of invariant manifolds for a wide class of `infinite dimensional' reaction-diffusion \spde{}s with linearly multiplicative noise; 
they built on earlier work by Bensoussan \& Flandoli~\cite{Bensoussan95} proving the existence and relevance of inertial manifolds of \spde{}s in a Hilbert space. Wang \& Duan~\cite{Wang06} also proved the existence of attractive slow manifolds for a wide class of \spde{}s but again only with linearly multiplicative noise.

Unfortunately, many interesting physical problems do not have Lipschitz nonlinearity; for example, the nonlinear advection in \bspde\ involves the unbounded operator~$\partial/\partial x$. I expect future theoretical developments should rigorously support this application. However, in the interim, let us proceed via a type of shadowing argument~\cite{Roberts05c}. The rapid dissipation of high wavenumber modes in applications such as \bspde\ indicates that the dynamics are close to finite dimensional: for example, I showed that resolving just ten subgrid modes in a modified \bspde\ on a small domain was sufficient to give the coefficients of the evolution on a \ssm\ correct to five digits~\cite[\S4]{Roberts05c}. By modifying any spatial derivatives in an~\spdee, such as the nonlinear advection in \bspde, to have a high wavenumber cutoff (a low-pass filter), as done in Section~\ref{sec:quad}, the dynamics of
an \spdee\ is effectively that of a Lipschitz, finite dimensional system. The \scm\ theorems of Boxler~\cite{Boxler89} then rigorously apply.

Others also invoke such a cutoff.  For example, Da~Prato \& Zabczyk~\cite[p.265]{DaPrato96} define a so-called `mollifier'~$M_R(k)$ which is identically~$k$ for $\|k\|<R$ and which is zero for $\|k\|>2R$\,.   In essence, the `wavenumber cutoff' (low-pass filter) of the previous paragraph performs the same regularising role as such a mollifer.

\paragraph{Coarse grain stochastic lattice dynamics}
Instead of a field continuous in space~$x$, suppose the microscale quantities and noise are known on a microscale lattice indexed by~$i$ as in the spatially discrete system~\eqref{eq:bssde}: for definiteness say $x_i=ih/2$  so that the coarse grid~$X_j$ has twice the spacing of the fine grid~$x_i$.  Then the microscale field is not the $L$-periodic continuum field~$u(x,t)$ but the $2m$-periodic discrete~$u(i,t)$. Similarly the noise is defined only at microscale lattice points.  Following earlier deterministic exploration~\cite{Roberts08c}, consider the stochastic lattice dynamics~\spdee\ when the dissipative operator is the specific second central difference $\cL u(i,t)=(4/h^2)\left[u(i+1,t)-2u(i,t)+u(i-1,t)\right]$ but the nonlinearity~$f$ and noise~$\phi$ may be quite general.  Embed the lattice dynamics onto $m$~overlapping elements to form the system~\eqref{eq:spdei} at internal lattice points and with the same coupling conditions~\eqref{eq:ibc}: the $j$th~element consists of microscale lattice points $x_i\in\{X_j, X_j\pm\frac12h,X_j\pm h\}$.  Such a discrete system has the same set of piecewise constant equilibria~$\vecs E_0$.  About each of these equilibria the linearised dynamics are the deterministic~\eqref{eq:spdel}.  The spectrum on each element, due to the second central difference~$\cL$, is then simply $\{0,-8/h^2,-16/h^2\}$ corresponding to microscale lattice eigenmodes (wavelets perhaps) on each element of $v_0=(1,1,1,1,1)$, $v_1=(0,-1,0,1,0)$ and $v_2=(1,-1,1,-1,1)$ respectively---as for the deterministic case~\cite{Roberts08c}.  Consequently, the Existence Corollary~\ref{lem:exist} then guarantees that for certain classes of lattice nonlinearity~$f$ and noise~$\phi$ there exists a \ssm\ describing an in principle exact closure, such as the approximate~\eqref{eq:bssdm}, on the coarse macroscale grid values of the microscale stochastic dynamics.

It will not escape your notice that such sound mapping of stochastic dynamics from a fine grid to a grid a factor of two coarser, such as \eqref{eq:bssde}$\mapsto$\eqref{eq:bssdm}, may be iterated across all grid scales on an entire multigrid hierarchy.  This approach has the potential to explore stochastic microscale dynamics at any scale, and to strongly relate dynamics across any scales. Thus this approach contributes to the need identified by Dolbow et al.~\cite[p.4]{Dolbow04} for ``representing information transfer across levels of scale'', and the more recent call by Brown et al.~\cite[p.14]{Brown08} for ``adaptive multiscale discrete stochastic simulation methods that are justified by theory and which can automatically partition the system into components at different scales''.

\subsection{The stochastic slow manifold captures emergent dynamics}

The second key property of \ssm{}s is that the evolution on the \ssm\ does capture the long term dynamics of the original \spdee, apart from exponentially decaying transients. For example, all solutions of \bspde{} close enough to the origin are are exponentially quickly described by the discrete \sde{}s~\eqref{eq:lowg}. This strong theoretical support for the model holds at finite element size~$h$---it ensures an accurate closure for the macroscale discretisation.

Boxler's Relevance Theorem~7.1(i) \cite{Boxler89}, or \cite[Theorem~4, p.213]{Roberts05c}, and Wang \& Duan's Theorem~A \cite{Wang06} are rephrased in the following corollary (see also Theorem~2.3(iii) by Bensoussan \& Flandolfi~\cite{Bensoussan95}) that, by the earlier argument, applies here to assure us that the dynamics on the \ssm\ emerges exponentially quickly.
 
\begin{corollary}[emergence]
	\label{lem:rel} 
	For the conditions of Corollary~\ref{lem:exist}, there is a finite neighbourhood~$N$ of the equilibria,~$(\mathbb E_0,\vec 0)$, such that for each solution~$u_j(x,t)$ of the \spde~\eqref{eq:spdei} starting and remaining in the neighbourhood~$N$,  there exists a solution~$\vec U(t)$ on the \ssm\ such that $\|u_j(x,t)-w_j(\vec U,x,t,\pars)\|\to0$ as $t\to\infty$ almost surely.
\end{corollary}

Boxler also assures us that the rate of decay to the \ssm\ is of the same magnitude as the decay of the gravest subgrid microscale mode~\cite[Theorem~7.1(i)]{Boxler89}, here $\exp(-\beta t)$\,. In \bspde\ for  example, on times significantly larger than a cross element diffusion time~$h^2/\pi^2$, the exponential transients decay and the \ssm\ model~\eqref{eq:lowg} describes the dynamics.  SImilarly, the transients of the microscale dynamics~\eqref{eq:bssde} decay on a time scale of~$h^2/8$ to the \ssm\ model~\eqref{eq:bssdm}.  For stochastic systems, such decay to the \ssm\ is most clearly seen via stochastic, normal form, coordinate transforms~\cite[e.g.]{Arnold98, Roberts06k, Roberts07d}.

Subtleties in this Emergence Corollary mislead some researchers, even when modelling deterministic systems. For example, Givon et al.~\cite{Givon04} discuss finite dimensional deterministic systems which linearly separate into slow modes~$x(t)$ and fast, stable modes~$y(t)$. They identify the existence of a slow invariant manifold $y=\eta(x)$ and the low~dimensional evolution on the manifold in the form $\dot X=L_1X+f(X,\eta(X))$\,. However, they \emph{assume}~\cite[p.R67, bottom]{Givon04} that the initial condition for the evolution on the slow manifold is simply $X(0)=x(0)$, and then consequently, just after their~(4.5), have to place undue restrictions on the possible initial conditions. However, the source of such restriction is that in general $X(0)\neq x(0)$. Physicists sometimes call the difference between $x(0)$~and~$X(0)$ the `initial slip'~\cite[e.g.]{Grad63}.  For stochastic systems, the correct nontrivial projection of initial conditions onto the \ssm\ may be constructed via stochastic normal forms~\cite[e.g.]{Arnold98, Roberts06k, Roberts07d} and applies in some finite domain around the subspace~$(\vecs E_0,\vec 0)$.

There are two caveats to our application of Corollary~\ref{lem:rel} to discretising \spde{}s. Firstly, although our constructed asymptotic series are global in the $m$~coarse macroscale amplitudes~$\vec U$, they are local in the parameters $\pars=(\sigma,\gamma,\alpha)$: the rigorous theoretical support applies in some finite neighbourhood of $\pars=\vec 0$\,. At this stage we have little information on the size of that neighbourhood. In particular, we evaluate the model when $\gamma=1$ to recover a discrete model for fully coupled elements; thus we require $\gamma=1$ to be in the finite neighbourhood of validity. Such validity has been demonstrated for the deterministic diffusion equation and Burgers' equation~\cite{Roberts98a}, but not yet for the stochastic case. Secondly, we cannot construct the \ssm\ and the evolution thereon exactly; it is difficult enough constructing asymptotic approximations such as the low order accuracy models \eqref{eq:diffm}, \eqref{eq:bssdm} and~\eqref{eq:lowg}. The models we develop and discuss have an error due to the finite truncation of the asymptotic approximations in the small parameters~$\pars$.

For example, the truncation in powers of the coupling parameter~$\gamma$ controls the width of the computational stencil for the discrete models. Due to the form of the coupling conditions~\eqref{eq:ibc}, nearest neighbour elements interactions are flagged by terms in~$\gamma^1$, whereas interactions with next to nearest neighbouring elements occur as $\gamma^2$~terms, and so on for higher powers. The low accuracy models \eqref{eq:diffm}, \eqref{eq:bssdm} and~\eqref{eq:lowg} are constructed with error~$\Ord{\gamma^2}$ and so only encapsulate interactions between the dynamics in an element and those of its two adjoining neighbours.  Construction to higher orders in coupling~$\gamma$ accounts for interactions between more neighbouring elements.

\section{Nonlinear dynamics have irreducible noise interactions}
\label{sec:quad}
Using the specific example of the forced \bspde, we now explore typical issues
arising in the macroscale discretisation of the quite general nonlinear \spdee, issues that arise for both differential and difference \spde{}s.

\subsection{Separate products of convolutions}
\label{sec:spc}

As detailed elsewhere~\cite[\S3.1]{Roberts06g}, in the iterative construction of the \ssm, $u_j=w_j(\vec U,x,t,\pars)=\uj+\cdots$ such that $d\uj/dt=g_j(\vec U,\pars)$\,, we use the residuals of the governing \spdee\ to drive corrections, indicated by primes, via 
\begin{equation}
    \D t{w'_j}-\cL_j{w'_j}+g'_jv_0
    =\text{residual}_{\eqref{eq:spde}}\,.
    \label{eq:diffd}
\end{equation}
In analysing a nonlinear \spdee, such as the stochastically forced \bspde, products of memory convolutions appear in the residual. Just seeking terms quadratic in the noise magnitude~$\sigma$ we generally have to deal with quadratic products of multiple convolutions in time.

\paragraph{Obtain corrections from residuals}

To cater for the general case, define multiple convolutions.  First, let the operator~$\Z{\frak k}$ denotes convolution over past history with $\exp[-\beta_{\frak k} t]$ where for brevity the fraktur~$\frak k$ denotes the pair~$jk$ corresponding to the $k$th~eigenvalue of the $j$th~element. That is,
\begin{equation}
        \Z{\frak k}\phi
	=\exp[-\beta_{\frak k}t]\star\phi(t) =\inti^t \exp[-\beta_{\frak k}(t-\tau)]
	\phi(\tau) \,d\tau\,;
    \label{eq:conv}
\end{equation}
recall that $\beta_{\frak k}$~is the (positive) decay rate of the $k$th~mode within the $j$th~element; $\beta_{jk}=k^2\pi^2/h^2$ for the \bspde, independent of the element.  Da~Prato \& Zabczyk~\cite[\S5.2.1]{DaPrato96} discuss the existence and continuity of such stochastic convolutions. Second, let $\Z{\vec {\frak k}}$ denote the operator of multiple convolutions in time where the subscript vector~$\vec {\frak k}$ indicates the decay rate of the corresponding convolution, that is, the operator
\begin{equation}
    \Z{\vec {\frak k}}=\Z{({\frak k}_1,{\frak k}_2,\ldots)} =\exp(-\beta_{{\frak k}_1}t)\star
    \exp(-\beta_{{\frak k}_2}t)\star \cdots\star
    \quad\text{and}\quad
    \Z{\,}=1\,,
    \label{eq:Z}
\end{equation}
in terms of the convolution~\eqref{eq:conv}; consequently
\begin{equation}
    \partial_t \Z{({\frak k}_1,{\frak k}_2,\ldots)}
    =-\beta_{{\frak k}_1} \Z{({\frak k}_1,{\frak k}_2,\ldots)}
    +\Z{({\frak k}_2,\ldots)}\,.
    \label{eq:Zt}
\end{equation}
The order of the convolutions does not matter~\cite{Roberts05c}; however, keeping intact the order of the convolutions seems useful to most easily cancel like terms in the residual of the governing \spde.

Earlier work~\cite[\S3.3--4]{Roberts06g} explored updates $w'_j$~and~$g'_j$ to the subgrid \ssm\ field and the discrete model \sde\ for deterministic terms~\cite[e.g.]{Roberts98a} and terms linear in the noise~$\phi$ in linear diffusion.   Summarising, in iteratively constructing the \ssm\ we encounter convolution integrals over the immediate past of the noise. Such fast time `memory' convolutions must be removed from the dynamics of the discretisation: Givon et al.~\cite[p.R59]{Givon04} similarly comment ``Memory. An important aim of any such algorithm is to choose~$P$ [the \ssm] in such a way that the dynamics in~$X$ [our grid values~$\vec U$] is memoryless.'' We simplify the discrete model tremendously by removing such `memory' convolutions as originally developed for \sde{}s by Coullet, Elphick \& Tirapegui~\cite{Coullet85}, Sri Namachchivaya \& Lin~\cite{Srinamachchivaya91}, and Roberts \& Chao~\cite{Chao95, Roberts03c}: the trick is to absorb the `memory' convolutions into the parametrisation of the \ssm, via~$w_j$, leaving the evolution on the \ssm,~$g_j$, to be convolution free. The results are discrete models, such as~\eqref{eq:diffm} and~\eqref{eq:bssdm}, where the noise is correlated across neighbouring elements despite the originally uncorrelated noise in the linear diffusion \spde.

Nonlinear spatially extended problems, such as the \bspde, need to adapt extra techniques from finite dimenisonal \sde{}s~\cite{Chao95, Roberts03c} and from \spde{}s on small spatial domains~\cite{Roberts05c}. For nonlinear problems we additionally have to solve for corrections quadratic in the noise for each term in the right-hand side of the general form
\begin{equation}
	\D t{w'_j}-\cL_j{w'_j}+g'_jv_0 = \sigma^2f_j(x)\Z{\vec
	{\frak \ell}}\phi_{j,n}\Z{\vec {\frak k}}\phi_{i,m}\,,
	\label{eq:wupnon}
\end{equation}
where $F_j(x)$ denotes complicated expressions encapsulating some of the influences of surrounding elements upon the subgrid structures within the $j$th~element.
Analogous to the treatment of linear stochastic terms~\cite[\S3.3--4]{Roberts06g}, two cases arise.
\begin{itemize}
    \item  Firstly, for each components of the subgrid
    structure~$F_j(x)$ in $\eigf_p(x-\xj)$ for wavenumber $p\geq
    1$\,, there is no difficulty in simply including in the
    correction to the subgrid field,~$w'_j$, the component
    \begin{displaymath}
        \sigma^2\eigf_p(x-\xj) \Z{jp}\left[ \Z{\vec
	\ell}\phi_{j,n}\Z{\vec{\frak k}}\phi_{i,m} \right]
    \end{displaymath}
    with its extra convolution in time.

	\item Secondly, for the component in~$F_j(x)$ that is constant across an element, the $\eigf_0(x-\xj)$~component, we separate the part of $\Z{\vec \ell}\phi_{j,n}\Z{\vec {\frak k}}\phi_{i,m}$ that has a bounded 	integral in time, and hence updates the subgrid field~$w'_j$, from the secular part that does not have a bounded integral and hence must update the model \sde\ through~$g'_j$.
\end{itemize}
In the considerations and convolutions for either case, the surrounding grid values that appear in the spatial structure forcing~$F_j$ are treated as constants as the time derivative in~\eqref{eq:wupnon} is the partial derivative of~$w_j(\vec U,x,t,\pars)$ keeping the grid values~$\vec U$ constant.  The $g'_j$~term in~\eqref{eq:wupnon} accounts for the time derivatives of grid values~$\vec U$.

\paragraph{Integrate by parts to separate}

Integration by parts introduced in earlier research~\cite{Roberts03c, Roberts05c} also here reduces all non-integrable convolutions to the canonical form of the convolution being entirely over one of the noises in a quadratic term.  Summarising, I define the canonical irreducible form to be $\phi_{j,n}\Z{\vec {\frak k}}\phi_{i,m}$.  Then, rewriting the convolution \ode~\eqref{eq:Zt} as $\beta_{\frak k}\Z{{\frak k}.\vec {\frak k}'}=-\partial_t{\Z{{\frak k}.\vec {\frak k}'}}+\Z{\vec {\frak k}'}$\,, where the vector of convolution parameters is decomposed as $\vec {\frak k}={\frak k}\cdot\vec {\frak k}'$ so that ${\frak k}$~is the first component of vector~$\vec {\frak k}$, and $\vec {\frak k}'$~is the vector (if any) of the second and subsequent components of vector~$\vec {\frak k}$, for any $\phi$~and~$\psi$, whether from the same element or not,
\begin{eqnarray*}
    \int \Z{{\frak k}\cdot\vec {\frak k}'}\phi \Z{\ell\cdot\vec \ell'}\psi\,dt
    &=&-\frac1{\beta_{\frak k}+\beta_\ell}\Z{{\frak k}\cdot\vec {\frak k}'}\phi  \Z{\ell\cdot\vec
    \ell'}\psi 
    \\&&{}
    +\frac1{\beta_{\frak k}+\beta_\ell}\int \Z{\vec {\frak k}'}\phi
    \Z{\ell\cdot\vec \ell'}\psi +\Z{{\frak k}\cdot\vec {\frak k}'}\phi \Z{\vec \ell'}\psi
    \,dt\,.
\end{eqnarray*}
Observe that each of the two components in the integrand on the right-hand side has one fewer convolutions than the initial integrand. Thus one repeats this integration by parts until terms of the cannonical form $\phi_{j,n}\Z{\vec {\frak k}}\phi_{i,m}$ in the integrand are reached. In this recursive process, assign all the integrated terms to update the subgrid stochastic field~$w'_j$. The irreducible terms remaining in the integrand, those in the form $\phi_{j,n}\Z{\vec {\frak k}}\phi_{i,m}$, must thus go to update the stochastic evolution~$g'_j$.

Computer algebra~\cite[\S6]{Roberts06b} readily implements these steps in the iteration to derive the asymptotic series of the \ssm\ of an \spdee.

\subsection{Odd noise highlights nonlinear noise interactions}
\label{sec:onhnni}

Subgrid modes with even wavenumber~$k$ are the only modes to affect the discrete model of linear diffusion~\cite[\S3.3]{Roberts06g}.  Consequently, restricting the forcing noise~$\phi(x,t)$ to have odd structure within each element highlights the \emph{nonlinear} dynamics of spatially extended dynamics of the \bspde, and clarifies the modelling of subgrid nonlinear stochastic effects in the general \spdee.

\begin{figure}
    \centering
    \includegraphics{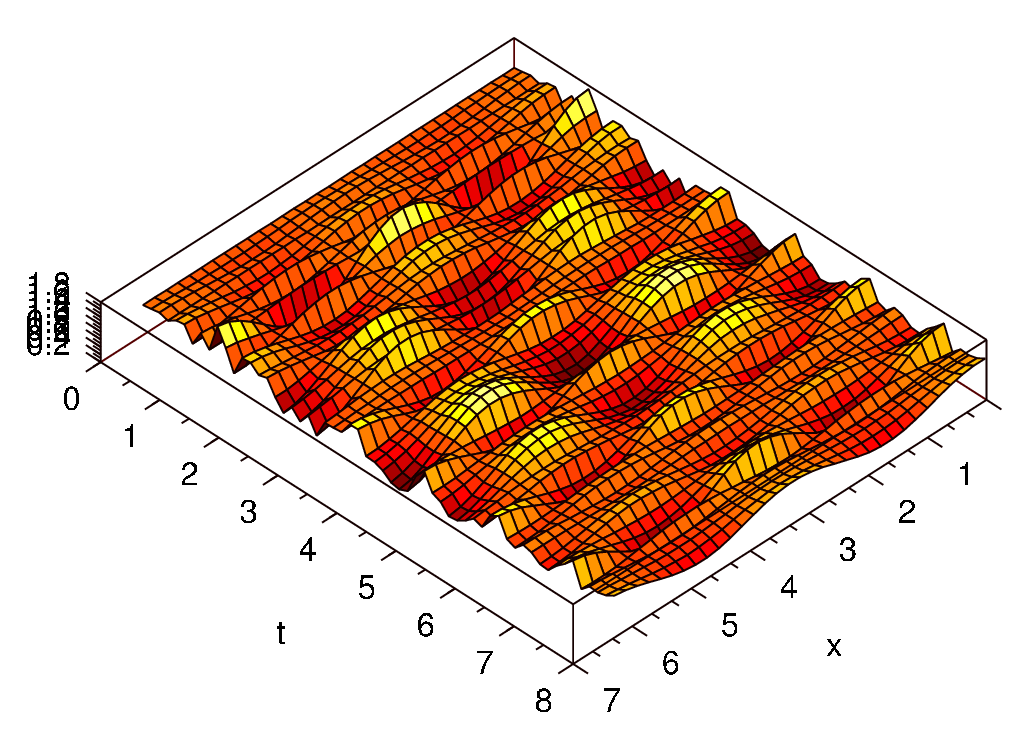} 
	\caption{microscale simulation of one realisation of the \bspde, $2\pi$-periodic and from initial condition $u(x,0)=1$\,, when the applied noise is simply $\sigma\psi(t)\cos 2x$ for $\sigma=1$ and nonlinearity $\alpha=\rat12$\,.  The 	macroscale grid points are at the nodes of this forcing, 	$x=\rat14\pi,\rat34\pi,\rat54\pi,\rat74\pi$\,, where the field~$u$ is relatively quiescent.}
    \label{fig:firsts}
\end{figure}

\begin{figure}
    \centering
    \begin{tabular}{c@{}c}
        \raisebox{20ex}{$u$} &
        \includegraphics{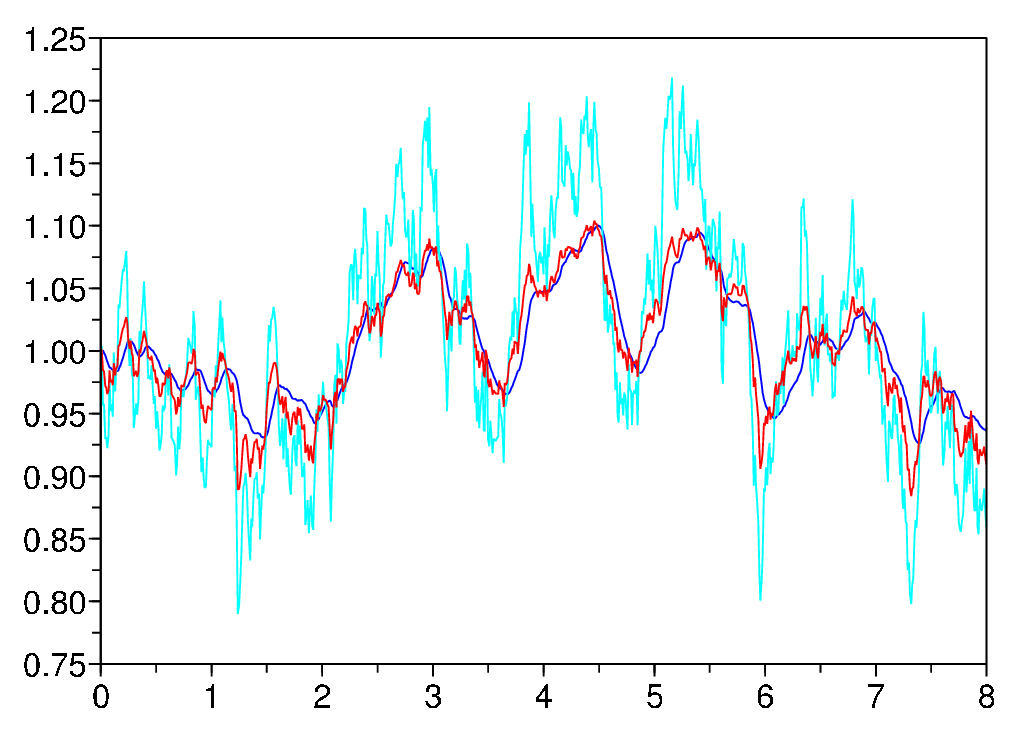} \\[-2ex]
        & time $t$
    \end{tabular}
	\caption{Compare the macroscale model~\eqref{eq:ssm1dt} with a microscale simulation for larger noise amplitude $\sigma=1$ and nonlinearity $\alpha=\rat12$ (as for Figure~\ref{fig:firsts}): 	blue, the microscale field $u(X_2,t)$ showing the subgrid noise carried by nonlinear advection past the grid point; cyan, the macroscale variable~$U_2(t)$; and red, the macroscale \ssm~\eqref{eq:ssm1g} at the grid point,~$u_2(X_2,t)=w_2(\vec U,X_2,t,\pars)$, for the shown~$U_2(t)$---the \ssm\ and the microscale simulation match well. }
    \label{fig:firstg2s1ah}
\end{figure}

\paragraph{One correlated noise}
The simplest nontrivial case of \bspde\ is when the noise~$\phi(x,t)$ has just the one odd Fourier component~$\sin\theta$ in each element, \emph{and} is perfectly but oppositely correlated in neighbouring elements.  That is, in this section set the noise in each element to be
\begin{equation}
    \phi(x,t)=(-1)^j\phi_1(t)\sin\theta\,,
    \label{eq:first}
\end{equation}
for the one `white noise'~$\phi_1(t)$. For example, Figure~\ref{fig:firsts} shows a microscale simulation of $2\pi$-periodic \bspde: as the macroscale element size $h=\pi/2$\,, the space-time structure of the noise~\eqref{eq:first} reduces to $\phi=\phi_1(t)\cos2x$. We proceed to model such microscale stochastic dynamics with four elements, with $h=\pi/2$\,, centred at $\xj=(j-\rat12)\pi/2$\,. Figure~\ref{fig:firsts} shows that the nodes of the stochastic forcing~\eqref{eq:first} are at these grid points~$\xj$.  Most methods for discretising \bspde\ on these elements, such as the simple collocation point sampling scheme~\eqref{eq:diffmm}, would sample or average over this noise and predict it has no influence whatsoever on the macroscale grid evolution. However, the subgrid scale nonlinear advection in \bspde\ carries the subgrid scale noise past the grid points and so generates fluctuations in the grid values: these fluctuations are hard to see in the simulation of Figure~\ref{fig:firsts}, but are clear in the plots of Figure~\ref{fig:firstg2s1ah}. In contrast to most methods for discretising \spde{}s, our holistic discretisation provides a systematic closure for the subgrid microscale dynamics and so predicts these fluctuations in the grid values.  

Computer algebra~\cite{Roberts06b} constructs the \ssm\ of \bspde\ with interelement coupling~\eqref{eq:ibc}: just modify the code to the specific forcing~\eqref{eq:first}.  In each element, and with noise induced, subgrid structures truncated to the first eight subgrid Fourier components~$1$, $\cosin\theta$, \ldots, $\cosin7\theta$\,, the subgrid microscale field is
\begin{eqnarray}&&
    u_j(x,t)=\uj
    +\gamma\Big[\frac{\theta}{\pi}\mud
    +\frac{\theta^2}{2\pi^2}\delta^2 \Big]\uj
    \pm\sigma\sin\theta\,\Z1\phi_1
    \nonumber\\&&{}
    \pm\sigma\alpha\uj\Big[
    \frac{2h}{\pi^2}\Z1
    -\frac4h\big( \rat13\cos2\theta\,\Z{2,1} -\rat1{15}\cos4\theta\,\Z{4,1}
    +\rat1{35}\cos6\theta\,\Z{6,1} \big)\Big]\phi_1
    \nonumber\\&&{}
    +\Ord{\sigma^3+\alpha^3+\gamma^{3/2}}\,,
    \label{eq:ssm1m}
\end{eqnarray}
where the upper alternative is for even~$j$, and the lower alternative for odd~$j$.  The terms in the second line of the \ssm~\eqref{eq:ssm1m} begin to account for the nonlinear advection and interaction of these subgrid spatial structures: these processes transform the additive forcing into multiplicative noise,~$\uj\phi_1$, with memory via the convolutions~$\Z{p,1}$.  Higher order terms in the coupling~$\gamma$, nonlinearity~$\alpha$ and noise magnitude~$\sigma$ are too onerous to record for the microscale subgrid structures.

A third \scm\ theorem of Boxler~\cite[Theorem~8.1]{Boxler89} justifies the asymptotic error reported in the \ssm~\eqref{eq:ssm1m}.   This error comes from the termination criterion of the iterative construction of the \ssm. Computer algebra~\cite{Roberts06b} iterates until the residual of the governing \bspde\ (or the fine lattice dynamics~\eqref{eq:bssde}), and the residual of the interelement coupling conditions~\eqref{eq:ibc}, are of some specified asymptotic order of error.  Boxler's~\cite{Boxler89} Theorem~8.1 then guarantees that the \ssm\ model constructed by the computer algebra has the same asymptotic order of error, as reported in the \ssm~\eqref{eq:ssm1m}.  

One may straightforwardly check by hand the error of an \ssm, such as~\eqref{eq:ssm1m}, by direct substitution into \bspde\ and the coupling conditions~\eqref{eq:ibc}.  However, the algebraic details are enormous, reflecting the intricate subgrid scale noise interactions and including the transformations outlined in Section~\ref{sec:spc}, and would fill many pages with otiose algebraic expressions.  Table~\ref{tbl:newt} tabulates the number of terms of various orders for the evolution on the \ssm\ of the \bspde: the number of terms describing the \ssm, higher order versions of~\eqref{eq:ssm1m}, is several orders of magnitude more.  Surely it is far better to leave such intricate detail to a computer, and focus instead on the stochastic model, its theoretical support, and its implications as we do here.

The construction of the \ssm~\eqref{eq:ssm1m} proceeds `hand-in-hand' with constructing the evolution on the \ssm.  Iterating to higher order in the small parameters~$\pars$, computer algebra~\cite{Roberts06b} gives that with the simple noise~\eqref{eq:first}, the emergent evolution of the grid values~$\uj(t)$ on the \ssm~\eqref{eq:ssm1m} is the system of \sde{}s
\begin{eqnarray}
    \duj &=& \frac\gamma{h^2}\delta^2\uj
    -\frac{\gamma^2}{12h^2}\delta^4\uj
    -\alpha\frac{\gamma}{h}\uj\mud\uj
    +\alpha^2\frac\gamma{12}\uj^2\delta^2\uj
    \nonumber\\&&{}
    \mp\sigma\alpha h\Big[ 
    \frac{2}{\pi^2} \uj
    +\gamma\big( .1028\, \uj +.0716\, \delta^2\uj \big) 
    -.00363\,\alpha^2h^2\uj^3 \Big]\phi_1
    \nonumber\\&&{}
    +\sigma^2\alpha^2\Big[ -\frac{8}{\pi^2}\uj\phi_1 \big(
    \rat1{15}\Z{2,1} +\rat1{255}\Z{4,1} +\rat1{1295}\Z{6,1} \big)\phi_1
    \nonumber\\&&\quad{}
    + 0.0195\, {h^2}\uj\phi_1 \Z1\phi_1 \Big]
    +\Ord{\sigma^5+\alpha^5+\gamma^{5/2}}\,.
    \label{eq:ssm1dt}
\end{eqnarray}
The terms in the first line form the deterministic holistic discretisation of Burgers' equation~\cite{Roberts98a}. The second line gives the terms linear in the noise: as for linear diffusion~\cite[\S3.2--3]{Roberts06g}, we remove all memory convolutions from the terms linear in noise. However, the terms in~\eqref{eq:ssm1dt} quadratic in noise,~$\sigma^2$, generally must contain memory convolutions (as discussed in Section~\ref{sec:spc}) in order to maintain the \sde{}s~\eqref{eq:ssm1dt} as a strong model of \bspde.

Figure~\ref{fig:firstg2s1ah} shows $U_2(t)$ for one simulation of the discrete model \sde~\eqref{eq:ssm1dt} for noise magnitudes~$\sigma=1$  (cyan curve). Compare these to the blue curve of the field~$u(X_2,t)$ at the corresponding grid point. Although the overall trends are roughly similar, the macroscale variable~$U_2(t)$ is markedly different to~$u(X_2,t)$. How then can the \sde~\eqref{eq:ssm1dt} be a strong model?  Resolve the difference by recalling that to eliminate memory convolutions we must abandon the freedom to impose precisely the meaning of the amplitudes~$\uj$ \cite[\S3, e.g.]{Roberts05c}. Thus, generally, $U_2(t)\neq u(X_2,t)$\,. Instead, the field at a grid point,~$u(\xj,t)$, is predicted by the
\ssm~\eqref{eq:ssm1m} evaluated at the grid points, namely
\begin{eqnarray}
    u_j(\xj,t)&=&\uj
    \pm\sigma\alpha\uj\Big[ \frac{2h}{\pi^2}\Z1 -\frac4h\big(
    \rat13\Z{2,1} -\rat1{15}\Z{4,1} +\rat1{35}\Z{6,1}
    \big)\Big]\phi_1
    \nonumber\\&&{}
    +\Ord{\sigma^4+\alpha^4+\gamma^2}\,.
    \label{eq:ssm1g}
\end{eqnarray}
Figure~\ref{fig:firstg2s1ah} plots the \ssm\ predicted grid value~\eqref{eq:ssm1g} in red and displays good agreement with the microscale simulation (blue). Evidently, \scm\ theory successfully supports discrete macroscale models of nonlinear \spde{}s.

\subsection{Strong models of stochastic dynamics are complicated}

Now restore independent and multiple noise processes in each element and consider the details of discretisations of the stochastically forced \bspde. Computer algebra~\cite{Roberts06b} derives the following leading terms in the asymptotic series of the model $d\uj/dt=g_j(\vec U,t,\pars)$. The large amount of algebraic detail reflects the enormous complexity of the multiple physical interactions acting on the subgrid microscale structures forced by the rich spectrum of stochastic noise. The arguments of the next section simplify the model significantly. Scan past the following model to the discussion.

Computer algebra derives that the element amplitudes~$\uj(t)$ evolve according to the system of \sde{}s
\begin{align}
    \duj={}&
    \gamma\frac1{h^2}\delta^2\uj 
    -\gamma^2\frac1{12h^2}\delta^4\uj
    -\gamma\alpha\frac1h\uj\mud\uj
    \nonumber\\&{}
    +\sigma\left\{ \vphantom{\frac83}
        \left[ 1 -\gamma\rat1{24}\delta^2
        +\gamma^2(\rat3{640}+\rat1{8\pi^4})\delta^4 \right]\phi_{j,0}
    \right.\nonumber\\&\left.\quad{}
        +\left[ \gamma\rat1{4\pi^2}\delta^2 -\gamma^2(\rat1{48\pi^2} 
        +\rat1{16\pi^4})\delta^4 \right]\phi_{j,2}
        -\alpha\frac{2h}{\pi^2}\uj\phi_{j,1}
    \right.\nonumber\\&\left.\quad{}
        +\alpha\gamma\frac{1}{h^2\pi^2}\left[  
            \uj\left( \rat8{\pi^2}\mud\phi_{j,0} 
                -\rat14\mud\phi_{j,2}
                +(\rat1{12}+\rat5{3\pi^2})\delta^2\phi_{j,1}
            \right)
    \right.\right.\nonumber\\&\left.\left.\qquad{}
            +\mud\uj\left( \rat14\phi_{j,2} 
                +(\rat16+\rat{10}{3\pi^2})\mud\phi_{j,1}
            \right)
    \right.\right.\nonumber\\&\left.\left.\qquad{}
            -\delta^2\uj\left( (\rat16+\rat1{3\pi^2})\phi_{j,1} 
                -(\rat1{24}+\rat5{6\pi^2})\delta^2\phi_{j,1} 
            \right)
        \right]
        -\alpha^2\frac{8h^2}{3\pi^4}\uj^2\phi_{j,0}
    \right\}
    \nonumber\\&{}
    +\sigma^2\left\{\vphantom{\frac11}
        \alpha\frac{h}{\pi^2}\left[ -2\phi_{j,0}\Z{1}\phi_{j,1}
            +\rat25\phi_{j,1}\Z{2}\phi_{j,2}
            +\rat25\phi_{j,2}\Z{1}\phi_{j,1} 
        \right]
    \right.\nonumber\\&\left.\quad{}
        +\alpha\gamma\frac1{h\pi^2}\left( -32\phi_{j,0}\Z{1,2}\mud
        -\rat45\phi_{j,1}\Z{2,2}\delta^2
        +\rat{32}5\phi_{j,2}\Z{1,2}\mud \right)\phi_{j,2}
    \right.\nonumber\\&\left.\quad{}
        +\alpha\gamma\frac{h}{\pi^2}\left[ 
            \phi_{j,0}\left( \rat8{\pi^2}\Z1\mud(\phi_{j,0}+\phi_{j,2}) 
            +(\rat1{12}+\rat5{3\pi^2})\Z1\delta^2\phi_{j,1}
    \right.\right.\right.\nonumber\\&\left.\left.\left.\quad\qquad{}
            -(\rat14+\rat8{\pi^2})\Z2\mud\phi_{j,2}
            \right)   
            +\phi_{j,1}\Z2\left( \rat15\delta^2\phi_{j,0}
            -(\rat1{20}+\rat{13}{150\pi^2})\phi_{j,2} \right)
    \right.\right.\nonumber\\&\left.\left.\qquad{}
            +\phi_{j,2}\left( -\rat8{5\pi^2}\Z1\mud(\phi_{j,0}+\phi_{j,2}) 
            -(\rat1{60}+\rat{17}{75\pi^2})\Z1\delta^2\phi_{j,1}
    \right.\right.\right.\nonumber\\&\left.\left.\left.\quad\qquad{}
            +(\rat18+\rat4{5\pi^2})\Z2\mud\phi_{j,2}
            \right)
    \right.\right.\nonumber\\&\left.\left.\qquad{}
            +\delta^2\phi_{j,0}\,\Z1\left(
            -(\rat1{12}+\rat2{15\pi^2})
            +(\rat1{24}+\rat5{6\pi^2})\delta^2 \right)\phi_{j,1}
    \right.\right.\nonumber\\&\left.\left.\qquad{}
            -\delta^2\phi_{j,1}\,\Z2\left(
            (\rat1{60}+\rat{17}{75\pi^2})
            +(\rat1{120}+\rat{17}{150\pi^2})\delta^2 \right)\phi_{j,2}
    \right.\right.\nonumber\\&\left.\left.\qquad{}
            -\delta^2\phi_{j,2}\,\Z1\left(
            (\rat1{20}+\rat{44}{75\pi^2})
            +(\rat1{120}+\rat{17}{150\pi^2})\delta^2 \right)\phi_{j,1}
    \right.\right.\nonumber\\&\left.\left.\qquad{}
            +\mud\phi_{j,0}\left(
            (\rat16+\rat{10}{3\pi^2})\Z1\mud\phi_{j,1}
            +(\rat14-\rat8{5\pi^2})\Z2\phi_{j,2} \right)
    \right.\right.\nonumber\\&\left.\left.\qquad{}
           -\mud\phi_{j,1} (\rat1{30}+\rat{34}{75\pi^2}) \Z2\mud\phi_{j,2}
    \right.\right.\nonumber\\&\left.\left.\qquad{}
            +\mud\phi_{j,2}\left(
            -(\rat1{30}+\rat{34}{75\pi^2})\Z1\mud\phi_{j,1}
            +(\rat18-\rat4{5\pi^2})\Z2\phi_{j,2} \right)
        \right]
    \right.\nonumber\\&\left.\quad{}
        +\alpha^2\frac1{\pi^2}\uj\left[ -\rat{16}3\phi_{j,0}\left(
            2\Z{1,2} +\rat{h^2}{\pi^2}\Z{2} \right)\phi_{j,2}
    \right.\right.\nonumber\\&\left.\left.\qquad{}
            -\rat8{15}\phi_{j,1}\left( \Z{2,1} -\rat{4h^2}{\pi^2}\Z{1}
            \right)\phi_{j,1} +\rat{16}{15}\phi_{j,2}\left( 2\Z{1,2}
            +\rat{h^2}{\pi^2}\Z{2} \right)\phi_{j,2}
        \right]
    \vphantom{\frac11}\right\}
    \nonumber\\&{}
    +\Ord{\sigma^3,\alpha^3+\gamma^3}\,.
    \label{eq:strongquad}
\end{align}

\paragraph{The model resolves noise, nonlinearity and inter-element interactions} The discrete model \sde~\eqref{eq:strongquad} is computed to residuals $\Ord{\sigma^3,\alpha^3+\gamma^3}$ and hence, supported by Boxler~\cite[Theorem~8.1]{Boxler89}, the model has the same order of error \cite[Theorem~5, p.213]{Roberts05c}. The \sde{}s~\eqref{eq:strongquad} build on earlier models of stochastic linear diffusion~\cite[\S3.4]{Roberts06g} recovered when one sets the nonlinearity parameter $\alpha=0$ in~\eqref{eq:strongquad}. The truncation to errors~$\Ord{\sigma^3}$ ensures the model retains the interesting mean effects generated by the quadratic noise interaction terms parametrised by~$\sigma^2$ seen in the last 14~lines of the \sde{}s~\eqref{eq:strongquad}. The truncation to error~$\Ord{\alpha^3+\gamma^3}$ resolves linear dynamics within and between next nearest neighbour elements, and nonlinear dynamics within and between nearest neighbour elements.

\begin{table}
    \newcommand{\non}[1]{\textit{#1}}
    \newcommand{\rrrr}{r@{\ \ }r@{\ \ }r@{\ \ }r}
    \centering
	\caption{number of terms in the evolution $d\uj/dt=g_j(\vec 	U,t,\pars)$ when only three Fourier modes are used for the subgrid stochastic structures: the numbers in \non{italics} report the terms evident in~\eqref{eq:strongquad}.  Expect many more terms when using more Fourier modes.  Blank entries are unknown.}
    \begin{tabular}{ccc}
    \begin{tabular}{r|\rrrr}
        \multicolumn{5}{c}{$\sigma^0$} \\[1ex]
        $\alpha^3$ & 0 & 0 & & \\
        $\alpha^2$ & \non0 & 3 & 14 &\\
        $\alpha^1$ & \non0 & \non2 & 8 & 19 \\
        $\alpha^0$ & \non0 & \non3 & \non5 & 7 \\
        \hline
        & $\gamma^0$ & $\gamma^1$ & $\gamma^2$ & $\gamma^3$ 
    \end{tabular}
    &
    \begin{tabular}{r|\rrrr}
        \multicolumn{5}{c}{$\sigma^1$} \\[1ex]
        $\alpha^3$ & 1 & 13 & & \\
        $\alpha^2$ & \non1 & 16 & 82 &\\
        $\alpha^1$ & \non1 & \non{11} & 45 & 93 \\
        $\alpha^0$ & \non1 & \non6 & \non{10} & 14 \\
        \hline
        & $\gamma^0$ & $\gamma^1$ & $\gamma^2$ & $\gamma^3$ 
    \end{tabular}
    &
    \begin{tabular}{r|\rrrr}
        \multicolumn{5}{c}{$\sigma^2$} \\[1ex]
        $\alpha^3$ & 9 & & & \\
        $\alpha^2$ & \non6 & 156 & &\\
        $\alpha^1$ & \non3 & \non{42} & 238 & \\
        $\alpha^0$ & \non0 & \non0 & \non0 & 0 \\
        \hline
        & $\gamma^0$ & $\gamma^1$ & $\gamma^2$ & $\gamma^3$ 
    \end{tabular}
    \end{tabular}
    \label{tbl:newt}
\end{table}

Computer memory~\cite{Roberts06b} currently limits us to the first few subgrid Fourier modes, $1,\sin\theta,\cos2\theta$, for this modelling of \bspde.  In modelling the microscale lattice dynamics~\eqref{eq:bssde} these three noise modes are complete, so in the application to coarse graining to~\eqref{eq:bssdm} the analysis accounts for all possible microscale noise. Table~\ref{tbl:newt} indicates the level of complexity of the multi-parameter asymptotic series via a type of Newton diagram.  The table reports the number of terms in various parts of the discrete model \sde\ $\dot \uj=g_j(\vec U,t,\pars)$ (there are vastly more terms describing the subgrid microscale structure~$w_j(\vec U,x,t,\pars)$ of the \ssm\ within each element). \emph{If we pursue either higher order truncations or more Fourier modes, then the complexity of the model increases alarmingly}.  Rational resolution of the subgrid scale stochastic interactions, in order to determine their macroscale effects,  suffers from a combinatorial explosion in terms. Thus, for the moment, truncate the model as in~\eqref{eq:strongquad}.

\paragraph{Abandon strong stochastic modelling}
The undesirable feature of the discrete model \sde{}s~\eqref{eq:strongquad} is the inescapable appearance in the quadratic noise terms of fast time convolutions, such as $\Z1 \phi_{j,1} =\exp(-\beta_1t)\star \phi_1$ and $\Z{1,2} \phi_{j,2} = \exp(-\beta_1t)\star \exp(-\beta_2t)\star \phi_{j,2}$\,. These require resolution of the subgrid fast time scales in order to maintain fidelity with the original \bspde{} and so would require incongruously small time steps for a supposedly slowly evolving model.  However, maintaining strong fidelity with the details of the white noise source~$\phi(x,t)$ is a pyrrhic victory when we are only interested in the relatively slow long term dynamics of the element amplitudes~$\uj(t)$. Instead, we need only those parts of the quadratic noise factors, such as $\phi_{j,0}\Z1 \phi_{j,1}$ and $\phi_{j,0}\Z{1,2} \phi_{j,2}$, that \emph{over the long macroscopic time scales} emerge as a mean drift and as new noise.  The next section develops such stochastic weak models in the context of a coarse grid discretisation.

\section{Stochastic resonance influences deterministic dynamics} 
\label{sec:sto}

Chao and Roberts~\cite{Chao95, Roberts03c, Roberts05c} argued that quadratic terms involving memory integrals of the noise were effectively new drift and new noise terms when viewed over long time scales  (as also noted by Drolet \& Vinal~\cite{Drolet97}). The arguments rely upon the noise being stochastic white noise. The strong model \sde~\eqref{eq:strongquad} faithfully tracks any given realisation of the original \bspde~\cite[Theorem~7.1(i), e.g.]{Boxler89} whether the forcing is deterministic or stochastic; however, this section proceeds to discuss how to derive a system of \sde{}s that weakly model on a macroscale a given spatially extended \spde. 

Analogously, Just et al.~\cite{Just01} argued that fast time deterministic chaos appears as noise when viewed over long time scales.

\subsection{Transform quadratic noise interactions}

In the strong model~\eqref{eq:strongquad} we need to understand and summarise the long term effects of the quadratic noises that appear in the form $\phi_j \Z{k} \phi_i$ and $\phi_j \Z{k,\ell} \phi_i$\,, where here $\phi_i$~and~$\phi_j$ represent the various possibilities for the components~$\phi_{j,k}$. The noises  may be independent ($s=0$) or they may be the same process ($s=1$) depending upon the term under consideration. Thus we explore the long term dynamics of Stratonovich stochastic processes $y_1$~and~$y_2$ defined via the \sde{}s
\begin{equation}
    \frac{dy_1}{dt}=\phi_j\Z{k}\phi_i
    \qtq{and}
    \frac{dy_2}{dt}=\phi_j\Z{\ell,k}\phi_i\,.
    \label{eq:stoy}
\end{equation}
Summarising earlier research~\cite{Chao95, Roberts03c, Roberts05c}, first name the two convolutions that appear in the nonlinear terms~\eqref{eq:stoy} as $z_1=\Z{k}\phi_i$ and $z_2=\Z{\ell,k}\phi_i$: then we must understand the long term properties of $y_1$~and~$y_2$ governed by the coupled system of \sde{}s
\begin{eqnarray}
    \dot y_1=z_1\phi_j\,,&&
    \dot z_1=-\beta_k z_1 +\phi_i\,,\nonumber\\
    \dot y_2=z_2\phi_j\,,&&
    \dot z_2=-\beta_\ell z_2 +z_1\,.
    \label{eq:bin}
\end{eqnarray}
Deterministic centre manifold theory applied to the Fokker--Planck equation for the system~\eqref{eq:bin} proves~\cite[\S4]{Roberts05c} that as time $t\to\infty$ the probability density function for the \sde~\eqref{eq:bin} tends \emph{exponentially quickly} to a quasi-stationary distribution~\cite[e.g.]{Pollett90}. The quasi-stationary distribution evolves according to a Kramers--Moyal equation which we interpret as approximating a Fokker--Planck equation for a system of \sde{}s (the neglected terms represent algebraically decaying non-Markovian effects among the $\vec y$~variables~\cite[equation~(11)]{Just01}). This established analysis~\cite[\S4]{Roberts05c} of the Fokker--Planck equation for system~\eqref{eq:stoy} models the system's long-time dynamics by the \sde{}s
\begin{equation}
    \frac{dy_1}{dt}=\rat12 s+\frac{\psi_1(t)}{\sqrt{2\beta_k}}
    \quad\mbox{and}\quad
    \frac{dy_2}{dt}=\frac{1}{\beta_k+\beta_\ell }\left(
    \frac{\psi_1(t)}{\sqrt{2\beta_k}}
    +\frac{\psi_2(t)}{\sqrt{2\beta_\ell }} \right).
    \label{eq:oosnn}
\end{equation}
As proved previously~\cite[Appendix~B]{Roberts05c}, and analogous to the argument of Just et al.~\cite[equation~(11)]{Just01}, the two $\psi_i(t)$ are \emph{new noises independent of $\phi_i$~and~$\phi_j$ over long time scales}. For the case of identical $\phi_i$~and~$\phi_j$ ($s=1$) there is a mean drift~$\rat12$ in the stochastic process~$y_1$; there is no mean drift in the other case of independent $\phi_i$~and~$\phi_j$ ($s=0$).

\subsection{Transform the complicated strong model to be usefully weak.} 
\label{sec:stow}

The challenge for our macroscale discretisation of a \spde\ is to handle the enourmous number of interacting noise processes in spatially extended nonlinear \spde{}s. For \bspde\ the decay rates $\beta_k=\pi^2k^2/h^2$. Thus, via various instances of the \sde{}s~\eqref{eq:oosnn}, to obtain a model for \emph{long time scales} we replace the quadratic noises in~\eqref{eq:strongquad} as follows:\footnote{Note that $\delta_{ij}$~and~$\delta_{mn}$, with its pair of subscripts, do \emph{not} denote a centred difference but rather denote the Dirac delta to cater for the self interaction of a noise when there is a mean drift effect ($s=1$), or not ($s=0$), as appropriate.}
\begin{eqnarray}
    \phi_{j,n}\Z{k}\phi_{i,m} &\mapsto& 
    \frac12\delta_{ij}\delta_{mn} +
    \frac{h}{k\pi\sqrt2}\psi_{nmk}(t)
    \,,\nonumber\\
    \phi_{j,n}\Z{k,\ell}\phi_{i,m} &\mapsto&
    \frac{h^3}{\pi^3(k^2+\ell^2)} \left[
    \frac1{\ell\sqrt2}\psi_{nm\ell}(t) 
    +\frac1{k\sqrt2}\psi_{nm\ell k}(t)
    \right]
    \,,\label{eq:xform}
\end{eqnarray}
where $\psi_{nmk}$~and~$\psi_{nm\ell k}$ are the effectively \emph{new
independent} white noises; that is, they are derivatives of new independent Wiener processes. I omit the subscripts of $i$~and~$j$ on~$\psi$ and henceforth on~$\phi$ because they are redundant when we record the model using centred mean and difference operators.

Computer algebra~\cite[\S9]{Roberts06b} implements the transformations~\eqref{eq:xform} to the strong model \sde{}s~\eqref{eq:strongquad} to derive the corresponding weak system of \sde{}s
\begin{align}
    \duj={}&
    \gamma\frac1{h^2}\delta^2\uj 
    -\gamma^2\frac1{12h^2}\delta^4\uj
    -\gamma\alpha\frac1h\uj\mud\uj
    \nonumber\\&{}
    +\sigma\left[ \phi_0
     - .2026\alpha h\uj \phi_1
     - .02738\alpha^2 h^2\uj^2\phi_2 \right]
    \nonumber\\&{}
    +\gamma\sigma\delta^2\left( - .04167\phi_0 
    + .02533\phi_2 \right)
    +\gamma^2\sigma\delta^4\left( .005971\phi_0 
        - .002752\phi_2 \right)
    \nonumber\\&{}
    +\alpha h\gamma\sigma\left\{ 
        \mud\uj\left[.02533\phi_2 + .05111\mud\phi_1 \right]
    \right.\nonumber\\&\left.\quad{}
        +\delta^2\uj\left[- .02031\phi_1 + .01278\delta^2\phi_1 \right]
    \right.\nonumber\\&\left.\quad{}
        +\uj\left[\mud\left( .08213\phi_0 - .02533\phi_2 \right) +
        .02555\delta^2\phi_1 \right]
    \right\}
    \nonumber\\&{}
    +\alpha h^2\sigma^2 \left( - .04561\psi_{011} 
        + .004561\psi_{122}
        + .009122\psi_{211} \right)
    \nonumber\\&{}
    +\alpha h^2\gamma\sigma^2\left[ 
        \mud_2\left( .01849\psi_{001} + .01849\psi_{021} -
        .01479\psi_{0212} -
        .01949\psi_{022} 
    \right.\right.\nonumber\\&\left.\left.\qquad{}
        - .003697\psi_{201} - .003697\psi_{221} +
        .002958\psi_{2212} + .003828\psi_{222} \right)
    \right.\nonumber\\&\left.\quad{}
        +\mud_1\left( .001002\psi_{022} + .0005011\psi_{222} \right)
    \right.\nonumber\\&\left.\quad{}
        +\mud_1\mud_2\left(.0115\psi_{011} - .0009038\psi_{122} -
        .001808\psi_{211} \right)
    \right.\nonumber\\&\left.\quad{}
        +\delta_2^2\left( .005752\psi_{011} + .0002311\psi_{102} -
        .0007858\psi_{122} 
    \right.\right.\nonumber\\&\left.\left.\qquad{}
        - .0001155\psi_{1222} - .0009038\psi_{211} \right)
    \right.\nonumber\\&\left.\quad{}
        +\delta_1^2\left( - .002209\psi_{011} - .0004519\psi_{122} -
        .002496\psi_{211} \right)
    \right.\nonumber\\&\left.\quad{}
        +\delta_1^2\delta_2^2\left( .002876\psi_{011} - .000226\psi_{122} -
        .0004519\psi_{211} \right)
   \right]
    \nonumber\\&{}
    +\alpha^2h^3\sigma^2\uj\left( - .004929\psi_{0212} 
        - .008626\psi_{022} + .004929\psi_{111} 
    \right.\nonumber\\&\left.\quad{}
        - .0002465\psi_{1112} - .0001232\psi_{112} + .0009859\psi_{2212} +
        .001725\psi_{222} \right)
    \nonumber\\&{}
    + .01643\alpha^2h^2\sigma^2\uj
    +\Ord{\sigma^3,\alpha^3+\gamma^3}\,,
    \label{eq:weakquad}
\end{align}
where $\mud_1$~and~$\delta_1^2$ denote differences in the first grid variable implicit in the noises~$\psi$, whereas $\mud_2$~and~$\delta_2^2$ denote differences in the second implicit grid variable. The value of the weak model \sde~\eqref{eq:weakquad} is that it has no fast time scale processes. The complexity of the weak model \sde~\eqref{eq:weakquad} reflects the intricacies of the inter-element interactions and the subgrid microscale processes resolved in this systematic approach to discretising \spde{}s.

The \sde{}s \eqref{eq:strongquad}~and~\eqref{eq:weakquad} may look very complicated, but the difference operators hide even more complexity in the stochastic components of the weak model \sde~\eqref{eq:weakquad}, and the strong model \sde~\eqref{eq:strongquad}. Any one apparent noise symbol~$\psi_{nmk}$ represents $5m$~independent noise sources over all the $m$~elements. In order to clarify all the discrete differences of~$\psi_{nmk}$ that appear, temporarily reinstate the implicit subscripts. For the $j$th~element we need the nine noise components $\psi_{jjnmk}$, $\psi_{j\pm1,jnmk}$, $\psi_{j,j\pm1,nmk}$~and~$\psi_{j\pm1,j\pm1,nmk}$ in order to compute all the differences that appear in the
\sde~\eqref{eq:weakquad}. Of these, seven of the noises are used in computing the differences in the $(j\pm1)$th~element, and two are also used in computing the differences in the $(j\pm2)$th~element. Consequently each of the $\sigma^2\gamma$~noises that appear in the
\sde~\eqref{eq:weakquad} actually represent, in nett effect, five independent noise sources for each element. Such terms reflect subtle cross-correlations between the stochastic dynamics within neighbouring finite elements.

\paragraph{Stochastic induced drift affects stability}
The terms quadratic in the noise magnitude, indicated by a factor~$\sigma^2$, are particularly complicated. With relatively small numerical coefficient, for many practical purposes we might ignore all except the last known term in the \sde~\eqref{eq:weakquad}, namely $+.01643\alpha^2h^2\sigma^2\uj$. The positive coefficient of this term shows that the self interactions of each of the many subgrid microscale noises act to promote mean growth of macroscale structures in the dynamics of \bspde.  Indeed, because of its potential importance, I included this mean effect in the introductory model \sde~\eqref{eq:lowg}.

Similarly, Boxler~\cite[p.544]{Boxler89}, Drolet \& Vinals~\cite{Drolet97, Drolet01},  Knobloch \& Weisenfeld~\cite{Knobloch83} and Vanden--Eijnden~\cite[p.68]{VandenEijnden05c} found stability modifying linear terms in their analyses of stochastically perturbed bifurcations and other nonlinear systems.

\subsection{Consolidate the emergent noises}

Orthonormalisation simplifies the representation of the effects of all the noise terms in the weak model \sde~\eqref{eq:weakquad}: this section reduces the 16~quadratic noises to just seven equivalent noise sources.  The novelty here is to do so for all the interacting elements. Because the noise terms appearing in~\eqref{eq:weakquad} are unknown in detail, we replace linearly independent combinations of them by one equivalent noise term, as for the discretisation of linear diffusion~\cite[\S3.4]{Roberts06g}. 

\begin{itemize}
\item Because of the severe truncation in the number of retained Fourier modes, there is no significant simplification possible in the terms linear in the noise magnitude~$\sigma$ that appear in~\eqref{eq:weakquad}.

\item Now turn to the quadratic noise terms in the
\sde~\eqref{eq:weakquad}. Computer algebra~\cite[\S9]{Roberts06b} extracts the eight different combinations of noises~$\psi$ in the \sde~\eqref{eq:weakquad}. Then a Gramm--Schmidt orthonormalisation of the vectors of coefficients is essentially a $QR$~decomposition of the transpose of the matrix of noise coefficients: namely, factor the noise contributions to~$R^TQ^T\vec\psi$ where $\vec\psi$~is the vector of noise processes, $Q^T$~is an orthogonal matrix, and $R^T$~is a lower triangular matrix. Then $\vec\chi=Q^T\vec\psi$ are a vector of new independent noise processes to replace~$\vec\psi$. For our weak model \sde~\eqref{eq:weakquad}, only the first seven rows of $R^T$ are non-zero, and hence only the first seven components of the new noises~$\vec\chi$ are significant. Thus seven new noises~$\vec\chi$, with coefficients in~$R^T$, replace the 16~noises~$\vec\psi$.
\end{itemize}
Computer algebra~\cite[\S9]{Roberts06b} computes the $QR$~factorisation of the quadratic noise coefficients in the weak model \sde~\eqref{eq:weakquad}, but extends to the case of four subgrid Fourier modes, instead of the three subgrid Fourier modes used to compute~\eqref{eq:weakquad}. There is a significant difference in the amount of detail: with truncation to four subgrid Fourier modes the weak model~\eqref{eq:weakquad} has 92~terms in its centred difference form; in comparison, with truncation to three Fourier modes the weak model~\eqref{eq:weakquad} has 53~terms in its centred difference form. However, upon replacing the quadratic noises~$\vec\psi$ by equivalent noises~$\vec\chi$ the resultant weak model has complexity largely independent of the number of retained subgrid Fourier modes. The model from four subgrid Fourier modes is
\begin{align}&
    \duj=
    \gamma\frac1{h^2}\delta^2\uj 
    -\gamma^2\frac1{12h^2}\delta^4\uj
    -\gamma\alpha\frac1h\uj\mud\uj
    \nonumber\\&{}
    +\sigma\left[ \phi_0
     +\alpha h\uj ( - .2026\phi_1 +.02252\phi_3)
     - .02555\alpha^2 h^2\uj^2\phi_2 \right]
    \nonumber\\&{}
    +\gamma\sigma\delta^2\left( - .04167\phi_0 
    + .02533\phi_2 \right)
    +\gamma^2\sigma\delta^4\left( .005971\phi_0 
        - .002752\phi_2 \right)
    \nonumber\\&{}
    +\alpha h\gamma\sigma\left\{ 
        \mud\uj\left[.02533\phi_2 + \mud(.05111\phi_1 -.00649\phi_3) \right]
    \right.\nonumber\\&\left.\quad{}
        +\delta^2\uj\left[(- .02031\phi_1- .0001769\phi_3) 
        +\delta^2(.01278\phi_1 -.001622\phi_3) \right]
    \right.\nonumber\\&\left.\quad{}
        +\uj\left[\mud\left( .08314\phi_0 - .02533\phi_2 \right) 
        +\delta^2( .02555\phi_1 -.003245\phi_3 ) \right]
    \right\}
    \nonumber\\&{}
    +.04681\alpha h^2\sigma^2 \chi_1
    \nonumber\\&{}
    +\alpha h^2\gamma\sigma^2\left[ 
        \mud_2\left( .02163\chi_2 +.02949\chi_3 \right)
    \right.\nonumber\\&\left.\quad{}
        +\mud_1\left( - .0006027\chi_2 - .000111\chi_3
        +.0008305\chi_4 \right)
    \right.\nonumber\\&\left.\quad{}
        +\mud_1\mud_2\left( - .01168\chi_1 +.000587\chi_5 \right)
    \right.\nonumber\\&\left.\quad{}
        +\delta_2^2\left( - .005875\chi_1 +.0001334\chi_5
        +.0004103\chi_6 \right)
    \right.\nonumber\\&\left.\quad{}
        +\delta_1^2\left( .001608\chi_1 - .002696\chi_5 -
        .0005192\chi_6 +.001116\chi_7 \right)
    \right.\nonumber\\&\left.\quad{}
        +\delta_1^2\delta_2^2\left( - .00292\chi_1 +.0001468\chi_5 \right)
   \right]
    \nonumber\\&{}
    +.01126\alpha^2h^3\sigma^2\uj\chi_2
    + .01751\alpha^2h^2\sigma^2\uj
    +\Ord{\sigma^3,\alpha^3+\gamma^3}\,.
    \label{eq:weakquadeq}
\end{align}
In this weak model, the new noises~$\chi_n$ \emph{implicitly} have two more subscripts to parametrise noise in pairs of nearby elements, as for~$\psi_{nmk}$: these reflect some of the subtle correlations between neighbouring elements.  The multitude of nonlinearity induced quadratic noise interactions have been replaced by just seven completely equivalent 	noise processes~$\vec\chi$ (although, as discussed earlier 	for~$\vec\psi$, these seven implicitly represent five times as many independent noise processes per element).  Again, because many of these quadratic noise terms have small coefficients, they may be neglected depending upon the parameter regime being explored, resulting in considerable simplification.  However, it is the multitude of such small noise interactions that lead to a significant mean drift effect.  The outcome of this section is that it is practicable to keep track of the multitude of interacting subgrid noises to extract important macroscale effects such as the mean drift term~$.01751 \alpha^2 h^2 \sigma^2 \uj$.

\section{Conclusion}
The virtue of the final weak model \sde~\eqref{eq:weakquadeq}, as also recognised by Just et al.~\cite{Just01}, is that we may accurately take large time steps as \emph{all} the fast dynamics are eliminated in the systematic closure. The critical innovation here is that we have demonstrated, via the particular example of \bspde, how it is feasible to analyse the net effect of many independent subgrid microscale stochastic effects, both within an element and between neighbouring elements, in a quite general \spdee. General formulae for modelling quadratic noise interactions~\cite{Roberts05c}, together with the iterative construction of stochastic slow manifold models~\cite{Roberts96a, Roberts06g}, empower us to model a wide range of \spde{}s, both differential and difference.

Theoretical support for the models comes from dividing the spatial domain into finite sized elements with coupling conditions~\eqref{eq:ibc}, invoking \scm\ theory~\cite{Boxler89}, and then systematically analysing the subgrid processes together with the appropriate physical coupling between the elements. This approach builds on success in discretely modelling deterministic \pde{}s~\cite[e.g.]{Roberts98a, Roberts00a, MacKenzie03}.  This theoretical support applies to stochastic difference equations on a spatial lattice, such as~\eqref{eq:bssde}, as well as to stochastic differential equations, such as \bspde.

What about spatial domains with physical boundary conditions at their extremes? The coupling parameter~$\gamma$ controls the information flow between adjacent elements; thus our truncation to a finite power in~$\gamma$ restricts the influence in the model of any physical boundaries to just a few elements near that physical boundary. The approach proposed here is based upon the \emph{local} dynamics on small elements while maintaining fidelity, via \scm\ theory, with the global dynamics of the original \spde. In the interior, the methods described here remain unchanged and thus produce identical model \sde{}s. The same methodology, but with different details will account for physical boundaries to produce a discrete model valid across the whole domain. Such modelling incorporating physical boundaries has already been shown for the deterministic Burgers' \pde~\cite{Roberts01b} and shear dispersion in a channel~\cite{MacKenzie03}.

Future research may find a useful simplification of the analysis used here if it can determine the mean drift terms, quadratic in~$\sigma^2$, without having to compute the other quadratic noise terms.  Techniques from stochastic averaging may help here~\cite[e.g.]{Pavliotis07}.

This approach to spatial discretisation of the \spdee\ may be extended to higher spatial dimensions as for deterministic \pde{}s~\cite{MacKenzie03, MacKenzie05, MacKenzie09b}.  Because of the need to decompose the stochastic residuals into eigenmodes on each element, the application to higher spatial dimensions are likely to require tessellating space into simple rectangular elements for \spde{}s.

\paragraph{Acknowledgement} The Australian Research Council supported this research with grants DP0774311 and DP0988738.

\bibliographystyle{unsrt}
\bibliography{ajr,bib}

\end{document}